\documentclass[10pt]{amsart}
\usepackage{url,graphicx}
\usepackage[dvips]{epsfig}
\usepackage{latexsym,color}
\usepackage{amscd,amssymb,verbatim}
\newcommand{\be}{\begin{enumerate}}
\newcommand{\ee}{\end{enumerate}}

\newcommand{\bd}{{\text{\rm Bd}\,}}

\newcommand{\cd}{{\mathcal D}}
\newcommand{\cf}{{\mathcal F}}

\newcommand{\cp}{{\mathcal P}}
\newcommand{\cs}{{\mathcal S}}
\newcommand{\rr}{{\mathbb R}}

\newcommand{\zz}{{\mathbb Z}_2}
\newcommand{\da}{\Delta}
\newcommand{\hra}{\hookrightarrow}
\newcommand{\lra}{\longrightarrow}
\newcommand{\bo}{\partial}

\newcommand{\dl}{{\textrm{dl}}}
\newcommand{\lk}{{\textrm{lk}}}
\newcommand{\idm}{{\textrm{id}}}
\newcommand{\aut}{{\textrm{Aut}\,}}

\newcommand{\la}{\leftarrow}
\newcommand{\ra}{\rightarrow}
\newcommand{\nin}{\noindent}
\newcommand{\pr}{\noindent{\bf Proof. }}
\newcommand{\sm}{\setminus}

\newcommand{\tmgn}{TM_{g,n}}

\newcommand{\tmgnb}{\tmgn^b}
\newcommand{\tmgne}{\tmgn^e}
\newcommand{\xgn}{X_{g,n}}
\newcommand{\cdgn}{\cd_{g,n}}

\newcommand{\colim}{{\textbf{colim}\,}}
\newcommand{\hocolim}{{\textbf{hocolim}\,}}
\newcommand{\hocolimm}{\textbf{\em{hocolim}}\,}
\newcommand{\fgraphs}{{\textbf{FGraphs}\,}}
\newcommand{\pgraphs}{{\textbf{PGraphs}\,}}
\newcommand{\ctop}{{\textbf{Top}\,}}
\newcommand{\ctopm}{{\textbf{\em{Top}}\,}}
\newcommand{\fixm}{\text{\em Fix}}
\newcommand{\fix}{\textrm{Fix}}

\newcommand{\ti}{\tilde}
\newcommand{\wti}{\widetilde}

\newtheorem{thm}{Theorem}[section]
\newtheorem{df}[thm]{Definition}
\newtheorem{lm}[thm]{Lemma}

\newtheorem{prop}[thm]{Proposition}
\newtheorem{conj}[thm]{Conjecture}
\newtheorem{rem}[thm]{Remark}

\numberwithin{equation}{section}
\numberwithin{figure}{section}
\numberwithin{table}{section}

\begin{document}

\title[Moduli spaces of tropical curves of higher genus]
{Moduli spaces of tropical curves of higher genus 
with marked points and homotopy colimits}

\author{Dmitry N. Kozlov}
\address{Department of Mathematics, University of Bremen, 28334 Bremen,
Federal Republic of Germany}
\email{dfk@math.uni-bremen.de}


\begin{abstract}
The main characters of this paper are the moduli spaces $\tmgn$ of
rational tropical curves of genus $g$ with $n$ marked points, with
$g\geq 2$.  We reduce the study of the homotopy type of these spaces
to the analysis of compact spaces $\xgn$, which in turn possess
natural representations as a~homotopy colimits of diagrams of
topological spaces over combinatorially defined generalized simplicial
complexes~$\da_g$, with the latter being interesting on their own
right.

We use these homotopy colimit representations to describe a~CW complex
decomposition for each $\xgn$. Furthermore, we use these developments,
coupled with some standard tools for working with homotopy colimits,
to perform an in-depth analysis of special cases of genus $2$ and $3$,
gaining a~complete understanding of the moduli spaces $X_{2,0}$,
$X_{2,1}$, $X_{2,2}$, and $X_{3,0}$, as well as a~partial
understanding of other cases, resulting in several open questions and
in further conjectures.
\end{abstract}

\maketitle

\section{Moduli spaces of tropical curves}

Tropical geometry is a~fairly recent new field within the broader
context of algebraic geometry. During the time which elapsed since its
inception, tropical geometry has already developed its language and
its methods, and has furthermore found numerous applications; we refer
the interested reader to \cite[Chapter 9]{cbms}, and more recently to
\cite{dfs, dy, mi2}, for both applications and the general background
information.

Some spaces arising in tropical geometry are of interest from the
point of view of algebraic topology as well. Often these have natural
definitions and fit well in more general structures. One such instance
is furnished by the moduli spaces of rational tropical curves of
genus~$g$ with $n$ marked points $\tmgn$, which were introduced by
Mikhalkin in~\cite{mi}, see also \cite{tms} for a~purely topological
definition.

Prior to this work, only the case of genus $1$ has been studied
systematically, see \cite{tms,tms2}, where, e.g., the homology groups
with coefficients in $\zz$ were computed for this family of spaces. In
this paper, we present an~in-depth analysis of the moduli spaces of
rational tropical curves of higher genus.

As a~first step, we complement the known shrinking bridges strong
deformation retraction, leading from $\tmgn$ to $\tmgnb$ as described
in~\cite[Section~3]{tms}, by a~further new simplification. On the
intuitive level, that newly discovered strong deformation retraction
increases the edge lengths proportionally, to reach the length of the
longest one, stopping when the edges which are strictly shorter than
the longest one form a~forest. This process is described in detail in
Section~\ref{sect:2}, where we reduce the space $\tmgnb$ to $\tmgne$
without changing the homotopy type (here $e$ stands for
``equalisation'').

By global scaling of edge lengths, we can replace $\tmgne$ by
a~compact space $\xgn$, and then proceed with finding more structure
in that new space. Hereafter, our main structural achievements are the
represention of $\xgn$ as a~homotopy colimit of a~diagram of
topological spaces in Section~\ref{sect:3}, and the derivation from
that representation of a~CW complex decomposition of $\xgn$ in
Section~\ref{sect:4}.

There are several byproducts of that development. The main one is
probably the discovery in Section~\ref{sect:3} of a~family of
generalized simplicial complexes $\da_g$, where $g$ is any natural
number.  The vertices of $\da_g$ are indexed by the isomorphism
classes of stable graphs of genus~$g$, and, more generally, the
simplices are indexed by filtrations by forests of these graphs.

Another byproduct is the introduction of cubical complexes (or, in
some terminology, of {\it generalized} cubical complexes) $C(G,\pi)$
in Section~\ref{sect:4}. Each such complex is defined using a~filtered
by forests stable graph as the input data. The interesting question of
connections between combinatorial properties of the input graphs and
geometry of the corresponding cubical complexes arises in a~natural
way.

We then use in Section~\ref{sect:5} and in Section~\ref{sect:6} the
previous developments to analyze the cases of genus $2$ and~$3$. In
the case when genus is equal to~$2$, we completely understand the
topology of moduli spaces when the number of marked points is $0$,
$1$, or~$2$. Furthermore, using an~Euler characteristic formula
derived using some standard enumeration under group action techniques,
we can show that these moduli spaces are almost never contractible.

Due to explosion in complexity, much less can be said in the case of
genus equal to~$3$. The main results here are the collapsibility of
the generalized simplicial complex $\da_3$, and the determination of
the asymptotics (with respect to the number of marked points $n$) of
the Euler characteristic of $\xgn$. The latter leads then to
a~conjecture concerning all genuses~$g$.

\section{An ``equalizing'' deformation retraction} \label{sect:2}

To start with, a few words on the terminology are in place. For
a~poset $P$, we denote the order complex by $\da(P)$. For a~graph $G$
we denote by $\da(G)$ the corresponding CW complex,
cf.\ \cite{tms}. For a~generalized simplicial complex $K$, we denote
by $\cf(K)$ its face poset. Recall that $\da(\cf(K))\cong\bd K$, see
e.g., \cite[p.\ 160, (10.4)]{CAT}.

We shall always work with finite undirected connected graphs only.
Furthermore, graphs satisfying the following additional properties
will play special role in this paper.
\begin{df}
A graph $G$ is called {\bf stable} if
\begin{itemize}
\item it has no bridges,
\item no vertex of $G$ has valency~$2$, unless it is adjacent to a~loop.
\end{itemize}
\end{df}

Due to the space constraints, we do not define the spaces $\tmgn$ and
$\tmgnb$ here, but rather refer to~\cite{mi} and~\cite{tms}.

\subsection{Filtered graphs} $\,$


\begin{df}
Given a~finite set $S$, an {\bf ordered set partition} of $S$ is
an~ordered tuple $\pi=(S_1,\dots,S_t)$, such that $S$ is a~disjoint
union of the subsets~$S_i$.
\end{df}

A~standard situation in which ordered set partitions arise is when we
have a~function $\varphi:S\ra\rr$, and we let $S_i$ be the non-empty
preimages $S_i:=\varphi^{-1}(x_i)$, $x_i\in\rr$, for $i=1,\dots,m$,
such that $x_1<x_2<\dots<x_m$.

\begin{df} \label{df:fgraph} $\,$

\nin (1) A {\bf filtered} graph is a~pair $(G,\pi)$, where $G$ is
a~graph and $\pi=(E_1,\dots,E_m)$ is an ordered set partition of
$E(G)$. We shall call $m$ the {\bf depth} of the filtration (or of the
filtered graph).

\nin (2) For such a filtered graph $(G,\pi)$ we shall say that $G$ is
     {\bf filtered by forests} if the subgraph induced by the edges
     $E_1\cup\dots\cup E_{m-1}$ is a~forest (i.e., contains no
     cycles). By convention, this condition is considered to be
     satisfied in the case $m=1$.
\end{df}

We remark that condition that $E_1\cup\dots\cup E_{m-1}$ induces
a~forest is equivalent to the condition that $E_1\cup\dots\cup E_j$
induces a~forest for all $j=0,\dots,m-1$.

Definition~\ref{df:fgraph} can be furthered by considering the
category of filtered graphs $\fgraphs$. To be specific, let $(G,\pi)$
and $(G',\pi')$ be two filtered graphs, say $\pi=(E_1,\dots,E_m)$ and
$\pi'=(E_1',\dots,E_{m'}')$. We have induced functions
$\rho:E(G)\ra[m]$ and $\rho':E(G')\ra[m']$ defined by saying that for
$e\in E(G)$ we have $e\in E_{\rho(e)}$, and for $e'\in E(G')$ we have
$e'\in E_{\rho'(e')}'$. 

\begin{df}
A~graph homomorphism $\varphi:G\ra G'$ is called a~{\bf filtered graph
  homomorphism} if for all $e_1,e_2\in E(G)$ we have the implication
\[\rho(e_1)\leq\rho(e_2)\Rightarrow\rho'(\varphi(e_1))\leq\rho'(\varphi(e_2)).\]
\end{df} 
\nin In other words, the function $\varphi$ preserves the (non-strict)
partial orders on $E(G)$ and $E(G')$ induced by the ordered set
partitions $\pi$ and~$\pi'$. 

Now, the objects of $\fgraphs$ are precisely all filtered graphs, and
morphisms are filtered graph homomorphisms. As usual in category
theory we call the invertible morphisms the {\it isomorphisms}. This
gives the notion of {\it isomorphic filtered graphs}. Being isomorphic
is an~equivalence relation, and thus we have a~notion of {\it
  isomorphism classes} of filtered graphs. Also, for every filtered
graph $(G,\pi)$ we get an~automorphism group $\aut(G,\pi)$, which
consists of all isomorphisms of $(G,\pi)$ with itself.

Before proceeding we would like to remark that omitting the order, and
considering pairs $(G,\pi)$, where $G$ is a~graph, and $\pi$ is
a~usual set partition of $E(G)$, will yield a~parallel concept of {\it
  partitioned graphs}. A {\it partitioned graph homomorphism} between
$(G,\pi)$ and $(G',\pi')$ is a~graph homomorphism $\varphi:G\ra G'$,
such that for any two edges $e_1,e_2\in E(G)$ the following condition
is satisfied: if $e_1$ and $e_2$ belong to the same block of $\pi$,
then $\varphi(e_1)$ and $\varphi(e_2)$ belong to the same block
of~$\pi'$. Accordingly, taking partitioned graphs as objects, and
partitioned graph homomorphisms as morphisms, we get a~category
$\pgraphs$. Furthermore, dropping the order of blocks in an~ordered
set partition induces a~forgetful functor from $\fgraphs$ to
$\pgraphs$.

To an~arbitrary metric graph $(G,l_G)$, the length function
$l_G:E(G)\ra(0,\infty)$ at hand allows us to associate a~filtered
graph $(G,\pi(G))$. This is done by considering an~ordered partition
$\pi(G)=(E_1,\dots,E_m)$ of the set $E(G)$ defined as follows:
\begin{itemize}
\item for every $1\leq i\leq m$, all the edges in $E_i$ have the same
  length, which we call $l(E_i)$;
\item we have $0<l(E_1)<\dots<l(E_m)$.
\end{itemize}
If this associated filtered graph is actually filtered by forests,
then we shall say that the original metric graph is filtered by
forests as well.

\begin{df}
The topological space $\tmgne$ is the subspace of $\tmgnb$ consisting
of the points whose representative metric graphs are filtered by
forests.
\end{df}


Our next goal is to show that as far as topology is concerned, it is
enough to consider the smaller space~$\tmgne$.

\subsection{The strong deformation retraction
  from $\tmgnb$ to $\tmgne$.} $\,$


\nin We start by describing an~explicit deformation
$\Phi:\tmgnb\times[0,1]\ra\tmgnb$. Let $x$ be a~point in $\tmgnb$,
represented by a~metric graph $(G,l_G)$, the corresponding ordered
partition $\pi(G)=(E_1,\dots,E_m)$, and the marking function
$p_G:[n]\ra\da(G)$. Let $k$ denote the maximal index, such that the
graph induced by $E_1\cup\dots\cup E_k$ is a~forest, i.e., has no
cycles.  If $E_1$ has a~cycle, then we set $k:=0$.

Informally speaking, the deformation $\Phi(x,-)$ should proceed as
follows:
\begin{itemize}
\item the lengths of the edges in $E_{k+1}\cup\dots\cup E_{m-1}$
  should increase so that the differences $l(E_m)-l(E_i)$, for
  $i=k+1,\dots,m-1$ decrease proportionally, and eventually the
  lengths $l(E_i)$ become equal to $l(E_m)$ at time $1$;
\item the lengths of the edges in $E_1\cup\dots\cup E_k$ should
increase  proportionally to the length $l(E_{k+1})$;
\item the positions of the points in the image of the marking function
  $p_G$ should also change proportionally to the increase of the
  lengths of the corresponding edges.
\end{itemize}

Formally, for $t\in[0,1]$, the point $\Phi(x,t)$ is given by the
representative $(G,l_{G}^t,p_{G}^t)$ which we now describe. To start
with, the graph itself (without the metric information taken into
account) is isomorphic to the graph~$G$. The length function $l_G^t$
is different, though for $t<1$ the corresponding ordered partition
$\pi(G)^t=(E_1,\dots,E_m)$ is the same, and we denote by $l_i^t$ the
value of $l_G^t$ on the edges in $E_i$; for $t=1$ the corresponding
ordered partition is $\pi(G)^1=(E_1,\dots,E_k,E_{k+1}\cup\dots\cup
E_m)$, and we also use the notation $l_i^1$, keeping in mind that
$l_{k+1}^1=\dots=l_m^1$.

More specifically, the values $l_i^t$ are given by the formulae:

\begin{equation} \label{eq:f1}
l_i^t:=l_i+t(l_m-l_i), \text{ for } k+1\leq i\leq m,
\end{equation}
and
\begin{equation} \label{eq:f2}
l_i^t:=l_i(1+t(l_m/l_{k+1}-1)), \text{ for } 1\leq i\leq k.
\end{equation}

Furthermore, if the point $p_G(i)$, for $i\in[n]$, is a~vertex of $G$,
then $p_G^t(i)=p_G(i)$, for all $t\in[0,1]$; else $p_G(i)$ is
an~internal point of some edge $e$ of $G$, in this case, the point
$p_G^t(i)$ belongs to the same edge, and has the same position {\it
  relative to} the length~$l_G^t(e)$.

\begin{thm}\label{thm:hmain}
The function $\Phi:\tmgnb\times[0,1]\ra\tmgnb$ defined above is
a~strong deformation retraction from $\tmgnb$ to $\tmgne$.
\end{thm}
\pr First, let the point $x\in\tmgnb$ be represented by a~metric graph
$(G,l_G)$, and let the corresponding ordered partition be
$\pi(G)=(E_1,\dots,E_m)$. When $k$ is the maximal index such that
$E_1\cup\dots\cup E_k$ has no cycles, it follows from \eqref{eq:f1}
that $l_{k+1}^1=\dots=l_m^1$, and hence
$\pi(G)^1=(E_1,\dots,E_k,E_{k+1}\cup\dots\cup E_m)$, in particular,
the point $\Phi(x,1)$ lies in $\tmgne$.

On the other hand, if that point $x$ lied in $\tmgne$ to start with,
then for the corresponding ordered partition $\pi(G)=(E_1,\dots,E_m)$
we would have that $E_1\cup\dots\cup E_{m-1}$ has no cycles, hence
$k=m-1$. It then follows from the formulae \eqref{eq:f1} and
\eqref{eq:f2} that $\Phi(x,t)=x$, for all $t\in[0,1]$.

It remains to show that the map $\Phi$ is continuous. Let a~point
$y\in\tmgne$ be represented by a~metric graph $(G,l_G)$. Choose
a~point $x\in\tmgnb$ and a~number $t\in[0,1]$, such that
$\Phi(x,t)=y$. We know that the point $x$ is represented by a~metric
graph $(G,l_G^x)$, and denote the corresponding ordered partition of
$E(G)$ by $\pi(G)=(E_1,\dots,E_m)$. Let $k$ as before denote the
maximal index such that the graph $E_1\cup\dots\cup E_k$ has no
cycles, and we set $l_i:=l(E_i)$, for $i=1,\dots,m$. 

Let $\varepsilon>0$ be an~arbitrary, sufficiently small number (say,
much smaller than all the non-zero distances between vertices and
marked points on $G$, in \cite{tms} we used the terminology ``the
admissible range''), and consider the $\varepsilon$-neighborhood
$N_\varepsilon(y)$. To show continuity of $\Phi$ we need to find
a~$\delta$-neighborhood
\[\wti N_\delta:=N_\delta(x)\times((t-\delta,t+\delta)\cap[0,1])\] 
of $(x,t)\in\tmgnb\times[0,1]$, such that $\Phi(\wti
N_\delta)\subseteq N_\varepsilon(x)$. We claim, that to do this, it is
enough to choose $\delta>0$ so that the inequality
\begin{equation} \label{eq:f3}
 \delta\cdot\frac{l_m}{l_{k+1}}<\varepsilon
\end{equation}
is satisfied. Note that the inequality does not depend on~$t$.

Indeed, pick a~point $\ti x\in N_\delta(x)$. By our definition of
$\tmgnb$ as a~topological space, which is given in detail in
\cite[Section~3]{tms}, it can be represented by the metric graph
$(\wti G,\ti l_G)$ together with the marking function $\ti p_G$, such
that
\begin{itemize}
\item the subgraph induced by the set $\Sigma$ of the edges whose
  length is less than $\delta$ has no cycles;
\item shrinking all the edges from $\Sigma$ inside $\ti G$ yields
  a~graph isomorphic to~$G$;
\item the lengths of all the other edges of $\ti G$ differ from
the lengths of the corresponding edges of $G$ by at most $\delta$;
\item the positions of the marked points, after the dilations of the
  edges they belong to, have also changed by at most~$\delta$.
\end{itemize}
It is now clear, that the point $\Phi(\ti x,\ti t)$ lies in the
$\varepsilon$-neighborhood of $y$, for any $\ti t\in[0,1]$. Indeed,
the point $\Phi(\ti x,\ti t)$ is represented by the graph isomorphic
to $\wti G$, with edges stretched by at most the factor $l_m/l_{k+1}$.
This shows that the edges of length less that $\delta$ are precisely
the ones which map to the edges of lengths less
than~$\varepsilon$. After they are shrunk we end up with a~graph which
is isomorphic to~$G$. Also, the lengths of other edges in the graph
representing $\Phi(\ti x,\ti t)$ differ from the lengths of the
corresponding edges in $G$ by at most $\delta l_m/l_{k+1}$, which by
our construction is less than~$\varepsilon$. Finally, the positions of
the marked points also did not change by more than $\varepsilon$.
Summarizing, we can conclude that $\Phi(\wti N_\delta)\subseteq
N_\varepsilon(x)$, and hence the function $\Phi$ is continuous. \qed

\vspace{5pt}

The space $\tmgne$ can be further simplified by normalizing the
lengths of the longest edges.

\begin{df} \label{df:xgn}
Let $\xgn$ denote the subset of $\tmgne$ consisting of points whose
representing metric graph $G$ has longest edges of length $1$. 
\end{df}

In our terminology above, if $\pi(G)=(E_1,\dots,E_m)$ is the ordered
partition corresponding to $G$, then the condition in
Definition~\ref{df:xgn} says that $l(E_m)=1$. Scaling metric graphs so
that their maximal edge length becomes equal to $1$ yields
homeomorphisms
\begin{equation}\label{eq:f4}
\tmgne\cong\xgn\times(0,\infty)\cong\xgn\times\rr,
\end{equation}
for all integers $g\geq 1$, $n\geq 0$. Since \eqref{eq:f4} implies
that the tropical moduli space $\tmgn$ is homotopy equivalent to
$\xgn$. From now on, we shall only work with the latter.

\begin{rem}
In~\cite{tms} and~\cite{tms2} the author defined and studied two
further variations of the space $\tmgn$. These are: $MG_{g,n}$ - the
moduli spaces of all metric graphs of genus $g$ with $n$ marked
points, and $MG_{g,n}^v$ - the subspace with the additional condition
that the marked points must be vertices. We notice here that
Theorem~\ref{thm:hmain} and its proof hold ad verbatim in these
generalized situations as well.
\end{rem}

\section{A stratification and homotopy colimit presentation} \label{sect:3}

In this section we shall describe how the tropical moduli space can be
replaced, up to homotopy equivalence, by a~manageable compact space,
which has a~nice presentation as a~homotopy colimit over a~certain
combinatorially defined generalized simplicial complex.

\subsection{The generalized simplicial complex of filtered by forests stable
graphs of genus~$g$} $\,$

\begin{df}
For an integer $g\geq 2$, let $\Sigma_g$ denote the set of all
isomorphism classes of filtered by forests stable graphs $(G,\pi)$,
where $G$ has genus~$g$. 
\end{df} 
Recall, that since the depth of isomorphic filtered graphs must be the
same, the notion of depth is well-defined for the isomorphism classes
as well.

\begin{df}
For an arbitrary integer $g\geq 2$ the generalized simplicial complex
$\da_g$ is defined as follows:

\nin \underline{\em the simplices:} the $m$-simplices of $\da_g$ are
indexed by the elements of $\Sigma_g$ of depth $m+1$;

\nin \underline{\em the boundary relation:} for an~$m$-simplex
$\sigma\in\da_g$, $m\geq 1$, let $(G,\pi)$ be a~representative of the
indexing element, say $\pi=(E_1,\dots, E_{m+1})$, then, the
representatives indexing the simplices on the boundary of $\sigma$ are
obtained by
\begin{itemize}
\item shrinking the edges from $E_1$ in $G$, and replacing $\pi$ with
  $(E_2,\dots,E_{m+1})$, or
\item keeping the graph $G$ intact, and merging two neighboring blocks
  $E_i$ and $E_{i+1}$, for $i=1,\dots,m$, in~$\pi$.
\end{itemize}
\end{df}

It is easy to see that the boundary relation in the generalized
simplicial complex $\da_g$ is well-defined by the description
above. Indeed, if $(G,\pi)$ is a~filtered by forests stable graph,
where $\pi=(E_1,\dots,E_{m+1})$, then
\begin{itemize}
\item merging blocks $E_i$ and $E_{i+1}$, for $i=1,\dots,m-1$, does
  not change the union of all blocks without the last one, hence the
  subgraph induced by that union remains being a~forest;
\item merging blocks $E_m$ and $E_{m+1}$ replaces the union
  $E_1\cup\dots\cup E_m$ with the union $E_1\cup\dots\cup E_{m-1}$,
  which is a~forest as well;
\item shrinking a~subposet inside a~forest still yields a~forest,
  hence $E_2\cup\dots\cup E_m$ induces a~forest in~$G/E_1$.
\end{itemize}

As a~special case, the vertices of $\da_g$ are indexed by the
isomorphism classes of filtered by forests stable graphs $(G,\pi)$,
where $\pi$ is a~$1$-tuple $(E(G))$, which unwinding all conditions,
simply translates to considering the isomorphism classes of stable
graphs of genus~$g$. The edges of $\da_g$ correspond to taking graphs
like that, choosing a~forest, and then identifying those forests which
map to each other under graph automorphisms. Vertices of every simplex
of $\da_g$ can be linearly ordered by the numbers of vertices in their
indexing stable graphs. This yields a~standard orientation on all the
simplices of $\da_g$, though we will not need that orientation in our
analysis.

For $g=1$ we make the convention that the generalized simplicial
complex $\da_1$ consists of a~single point, which corresponds to the
graph with one vertex and a~loop at that vertex. This is consistent
with our definition of a~stable graph.

For $g=2$ we have two isomorphism classes of stable graphs of
genus~$g$: a~graph consisting of one vertex with two loops attached,
and a~graph consisting of two vertices connected by three edges. In
the first case there are no non-empty forests. In the second case,
a~non-empty forest is given by taking one of the edges, and this
choice is unique up to graph automorphisms. Thus, the generalized
simplicial complex $\da_2$ is simply a~$1$-simplex, see
Figure~\ref{fig:g21}.

The case $g=3$ is more complicated, see Figures~\ref{fig:g32}
and~\ref{fig:g33}. It will be shown in Section~\ref{sect:g3} that the
complex $\da_3$ is collapsible.

\begin{figure}[hbt]
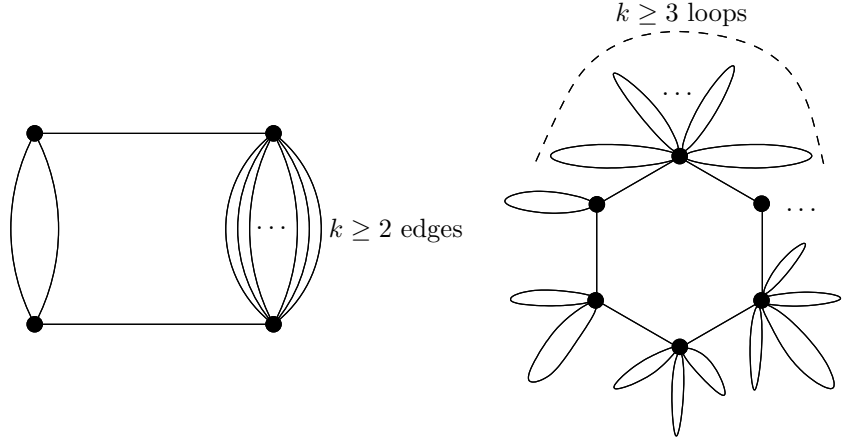

\begin{center}
  \begin{picture}(0,0)%
    \includegraphics{f1.pstex}%
  \end{picture}%
  \input{f1.pstex_t}%
  
\end{center}
\caption{On the left hand side we show a graph of genus $k+1$ which
  has $2$ non-isomorphic spanning forests, for any $k\geq 2$. On the
  right hand side, we show graph consisting of a $k$-gon, $k\geq 3$,
  with $i$ loops attached to the vertex number $i$. This graph has $k$
  non-isomorphic spanning forests. }
\label{fig:f1}
\end{figure}

Let us now look at some properties of generalized simplicial complex
$\da_g$ valid for all~$g$. First of all $\da_g$ connected. To see
this, let $v_l$ denote the vertex which is indexed by the graph $L$
with one vertex and $g$ loops. Let $w$ be some other vertex of $\da_g$
represented by a~stable graph~$H$. The edges connecting $w$ with $v_l$
correspond to various spanning trees in the graph $H$, with two trees
$T_1$ and $T_2$ corresponding to the same edge if and only if there
exists a~graph automorphism of $H$ which transforms $T_1$ to $T_2$.
Since at least one spanning tree always exists, there is at least one
edge connecting $w$ with $v_l$, in particular, the complex $\da_g$ is
connected. The example on the left hand side of Figure~\ref{fig:f1},
shows that $\da_g$ is a~simplicial complex only in cases $g=1,2$,
while the example on the right hand side of Figure~\ref{fig:f1} shows
that there could be arbitrary many edges between vertices for higher
genuses.

\begin{prop} \label{prop:3.4}
For arbitrary integer $g\geq 1$, the generalized simplicial
complex~$\da_g$ is pure. It has dimension $0$ for $g=1$, and $2g-3$
for $g\geq 2$.
\end{prop}

\pr The statement is true for $g=1$, so we assume that $g\geq
2$. Consider a~maximal simplex $\sigma$ of $\da_g$, say it is
represented by a~filtered by forests stable graph $(G,\pi)$, where
$\pi=(E_1,\dots,E_m)$. By our construction, the dimension of $\sigma$
is~$m-1$. Since $\sigma$ is maximal, and we are always allowed to
split the first $m-1$ blocks of $\pi$ into smaller ones, we may assume
that $|E_1|=\dots=|E_{m-1}|=1$. Furthermore, since we are allowed to
split of an edge from $E_m$ as long as its union with the other blocks
forms a~forest, we may assume that $E_1\cup\dots\cup E_{m-1}$ induces
a~spanning forest.

We now claim that all vertices of $G$ have valency $3$. Assume this is
not the case, and take a~vertex $w$ which has valency at least $4$.
We consider two cases. 

\nin {\bf Case 1.} There exists a loop adjacent to $w$. Let $l$ denote
that loop, and let $e_1,\dots,e_t$ denote the other edges adjacent to
$w$. Here, by our assumptions, we have $t\geq 2$, and loops (other
than $l$) will appear twice in that list. Let $G'$ be a new graph,
which is obtained from $G$ by 
\begin{itemize} 
\item replacing the vertex $w$ with new vertices $w_1$ and $w_2$;
\item connecting the edges $e_1,\dots,e_{t-1}$ to $w_1$, $e_t$ to
  $w_2$, and replacing the loop $l$ with two new edges $l_1$ and
  $l_2$, both connecting $w_1$ and $w_2$.
\end{itemize}

\begin{figure}[hbt]
\begin{center}
  \begin{picture}(0,0)%
    \includegraphics{f2.pstex}%
  \end{picture}%
  \input{f2.pstex_t}%
  
\end{center}
\caption{The graph transformation used in Case 1 of the proof of
  Proposition~\ref{prop:3.4}}
\label{fig:f2}
\end{figure}

This transformation is shown graphically on Figure~\ref{fig:f2}.  We
see that shrinking the edge $l_2$ in the graph $G'$ will yield a~graph
isomorphic to $G$, and that $G$ and $G'$ have the same
genus. Furthermore, $G'$ has no bridges and no vertices of valency
$2$, so it is a~stable graph. Setting $\pi':=(\{l_2\},E_1,\dots,E_m)$,
we get a~filtered stable graph $(G',\pi')$. The graph
$G'/(\{l_2\},E_1,\dots,E_{m-1})$ is isomorphic to
$G/(E_1,\dots,E_{m-1})$, hence has genus $g$, implying that
$\{l_2\},E_1,\dots,E_{m-1}$ induces a~forest. Summarizing, we conclude
that $(G',\pi')$ is a~filtered by forests stable graph which is
indexing a~simplex $\tau$, such that $\dim\tau=\dim\sigma+1$, and
$\sigma\subset\tau$, contradicting the fact that $\sigma$ is a~maximal
simplex.

\begin{figure}[hbt]
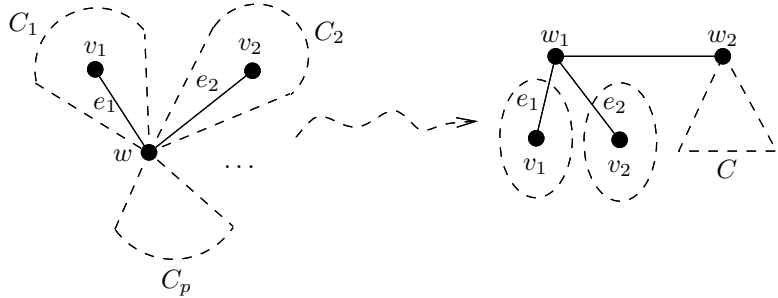

\begin{center}
  \begin{picture}(0,0)%
    \includegraphics{f4.pstex}%
  \end{picture}%
  \input{f4.pstex_t}%
  
\end{center}
\caption{The deformation of the graph $G$.}
\label{fig:f4}
\end{figure}

\nin {\bf Case 2.} All the edges adjacent to $w$ are not loops. Denote
these edges $e_1,\dots,e_t$, and assume that $e_i$ connects $w$ to
a~vertex $v_i$, for $i=1,\dots,t$.  Let $C_1,\dots,C_p$ denote the
connected components of the graph obtained from $G$ by removing the
vertex $w$ and all the adjacent edges. It is important to remark that,
for all $i=1,\dots,p$, the valency of $w$ in the subgraph of $G$
induced by $C_i\cup\{w\}$ is at least $2$, as otherwise the single
edge adjacent to $w$ in $C_i\cup\{w\}$ would have been a~bridge in the
graph~$G$. Let $e_1$ be an~edge adjacent to $w$ in $C_1\cup\{w\}$, and
let $e_2$ be an~edge adjacent to $w$ in $C_p\cup\{w\}$, see the left
hand side of Figure~\ref{fig:f4}; we might have $p=1$, but that does
not change the argument.  Let now $G'$ be a new graph, which is
obtained from $G$ by
\begin{itemize} 
\item replacing the vertex $w$ with new vertices $w_1$ and $w_2$;
\item connecting the edges $e_1$ and $e_2$ to $w_1$, $e_3,\dots,e_t$
  to $w_2$, and adding a~new edge $e$ connecting $w_1$ and $w_2$.
\end{itemize}

\begin{figure}[hbt]
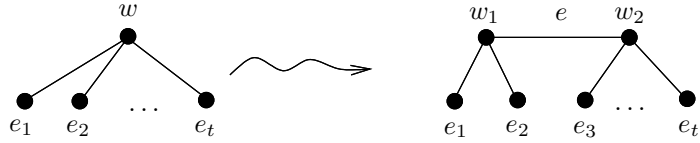

\begin{center}
  \begin{picture}(0,0)%
    \includegraphics{f3.pstex}%
  \end{picture}%
  \input{f3.pstex_t}%
  
\end{center}
\caption{The graph transformation used in Case 2 of the proof of
  Proposition~\ref{prop:3.4}}
\label{fig:f3}
\end{figure}

This transformation is shown on Figure~\ref{fig:f3}. Set
$\pi':=(\{e\},E_1,\dots,E_m)$. Essentially with the same proof as in
the first case, we see that $(G',\pi')$ is filtered by forests stable
graph, and that shrinking $e$ will yield $(G,\pi)$. The only fact
which needs to be verified is that $e$ is not a~bridge in $G'$.
Indeed, on one hand, every vertex in
$V(G')\sm\{w_1,w_2\}=V(G)\sm\{w\}$ is in the same connected component
of $G'\sm\{e\}$ as $w_2$, since it is in some $C_i$ and the valency of
$w$ in $C_i\cup\{w\}$ is at least $2$, hence at least one of these
edges is different from $e_1$, $e_2$ (when $p=1$, both $e_1$ and $e_2$
are adjacent to $w$, but then the valency of $w$ in $C_1\cup\{w\}$ is
at least $4$). On the other hand $w_1$ is connected to two other
vertices in $G'\sm\{e\}$, hence is also in the same connected
component of $G'\sm\{e\}$ as $w_2$. We conclude again that $(G',\pi')$
is a~filtered by forests stable graph which is indexing a~simplex
$\tau$, such that $\dim\tau=\dim\sigma+1$, and $\sigma\subset\tau$,
contradicting the fact that $\sigma$ is a~maximal simplex.

It is now easy to finish the proof. Assume $G$ has $v$ vertices and
$e$ edges. Then we have $g=e-v+1$. On the other hand, we have $2e=3v$.
It follows that $v=2g-2$ and $e=3g-3$. Since $E_1\cup\dots\cup
E_{m-1}$ induces a~spanning forest, we have $m-1=v-1$, and it follows
that $\dim\sigma=2g-3$.  \qed

\vspace{5pt}
 
We are now ready for the main result of this subsection.

\begin{thm} \label{thm:xg0}
For arbitrary integer $g\geq 1$, the generalized simplicial complex
$\da_g$ is homeomorphic to the space~$X_{g,0}$.
\end{thm}

\pr We shall define a~map $\rho:X_{g,0}\ra\da_g$. Let a~point $x\in
X_{g,0}$ be represented by a~filtered by forests stable metric graph
$(G,l_G)$, with the corresponding ordered partition
$\pi(G)=(E_1,\dots,E_m)$. The point $\rho(x)$ belongs to the simplex
of $\da_g$ indexed by the isomorphism class of $(G,\pi(G))$. Its
vertices $\{v_1,\dots,v_m\}$ are indexed by the graphs $G_i$, for
$i=1,\dots,m$, where $G_i$ is obtained from $G$ by shrinking all the
edges in the set $E_1\cup\dots\cup E_{i-1}$. The coordinate $d_i$ of
the vertex $v_i$ is given by $d_i:=l_i-l_{i-1}$, where $l_i=l(E_i)$,
$l_m=1$, and by convention $l_0:=0$. Clearly, we have $d_i>0$, for all
$i=1,\dots,m$, and $d_1+\dots+d_m=l(E_m)=1$.

Let us check that the described map $\rho:x\mapsto d_1 v_1+\dots+d_m
v_m$ is indeed a~homeomorphism $\rho:X_{g,0}\ra\da_g$. First, it is
clearly bijective as the $m$-tuple transformation
$(l_1,\dots,l_m)\mapsto(l_1,l_2-l_1,\dots,l_m-l_{m-1})$ has the
inverse $(d_1,\dots,d_m)\mapsto(d_1,d_1+d_2,\dots,d_1+\dots+d_m)$,
which defines $\rho^{-1}:\da_g\ra X_{g,0}$. Second, both $\rho$ and
its inverse are continuous. Indeed, for $x\in X_{g,0}$ represented by
$(G,l_G)$, with the corresponding ordered partition
$\pi(G)=(E_1,\dots,E_m)$, the representatives for the points in
a~small neighborhood are obtained by a~combination of the following 3
steps:
\begin{enumerate}
\item[(1)] changing edge lengths without splitting blocks: for some
  $1\leq i\leq m-1$, the edges in $E_i$ get the length
  $l_i+\varepsilon$ instead of $l_i$, where $\varepsilon$ is a~small,
  not necessarily positive number;
\item[(2)] splitting of a~block: the set $E_i$ gets split into
  a~disjoint union of sets $A$ and $B$, with edges in $A$ keeping the
  length $l_i$, and edges in $B$ getting the length
  $l_i+\varepsilon$, where $\varepsilon$ is a~small positive number;
\item[(3)] adding a~block of short edges: we replace some of the
  vertices of $G$ with trees, and add a~block $E$ consisting of all
  the new edges, with all the edges in $E$ getting the length
  $\varepsilon$, where $\varepsilon$ is again a~small positive number.
\end{enumerate}
On the side of $\da_g$ these 3 steps can be interpreted as follows.
Step (1) corresponds to the case when $v_1,\dots,v_m$ stay the same,
whereas coordinates $d_i$ and $d_{i+1}$ change to $d_i+\varepsilon$
and $d_{i+1}-\varepsilon$. Step (2) corresponds to the situation where
we add a new vertex $w$ indexed by the graph obtained from $G$ by
shrinking $E_1\cup\dots\cup E_{i-1}\cup A$; this vertex gets
coordinate $\varepsilon$, whereas the vertex $v_i$ gets coordinate
$d_i-\varepsilon$.  Finally, step (3) corresponds to the situation,
where we add a new vertex $w$ indexed by the graph obtained from $G$
by shrinking $E$; this vertex gets coordinate $\varepsilon$, whereas
the vertex $v_1$ gets coordinate $d_1-\varepsilon$.  \qed


\subsection{A diagram over the face poset of $\da_g$} $\,$

\nin For a finite set $S$, let $X_{g,S}$ denote the moduli space of
graphs, which is just like $\xgn$, but the labels for the marked
points are drawn from the set $S$, i.e., the labeling function is
$p_G:S\ra\da(G)$. In particular, $\xgn=X_{g,[n]}$. For any injective
set map $\iota:S\hra T$, we get an induced map $\ti\iota:X_{g,T}\ra
X_{g,S}$, which takes every metric graph $G$ to itself, but changes
the old labeling function $p_G:T\ra\da(G)$ to the new one
$q_G:S\ra\da(G)$ defined by $q_G(s):=p_G(\iota(s))$, for~$s\in S$.

Composition of two injective maps $\iota_1:S_1\hra S_2$ and
$\iota_2:S_2\hra S_3$ corresponds to the composition of the induced
maps $\ti\iota_1:X_{g,S_2}\ra X_{g,S_1}$ and $\ti\iota_2:X_{g,S_3}\ra
X_{g,S_2}$, as we have
$\wti{\iota_2\circ\iota_1}=\ti\iota_1\circ\ti\iota_2$. This means that
the map $\lambda:S\mapsto X_{g,S}$ yields a~contravariant functor from
the category {\bf Inj} of finite sets and injective set maps to the
category {\bf Top} of topological spaces and continuous maps. The
category {\bf Inj} has an~initial element $\emptyset$, hence, for
every finite set $S$ we get the unique induced map
$\lambda(\emptyset\hra S):X_{g,S}\ra X_{g,0}$.  Clearly, the induced
map $\lambda(\emptyset\hra S)$ is simply the forgetful map which
``forgets'' all the labels.

We have described above and proved in Theorem~\ref{thm:xg0} that the
space $X_{g,0}$ has a~natural structure of a~generalized simplicial
complex. This can also be viewed as a~stratification of $X_{g,0}$ by
open simplices. It is then a~standard construction to consider the
stratification on $\xgn$ induced by that forgetful map: the strata are
simply the preimages of the strata of $X_{g,0}$ under
$\lambda(\emptyset\hra S)$.

To formulate the second main result of this section we need to recall
some terminology of diagrams and homotopy colimits. We shall only
consider simplified version which we need here, and refer the reader
to \cite[Chapter 15]{CAT} for a~more complete coverage of the general
situation.

\begin{df} \label{df:diag}
 Let $P$ be a~poset. A~{\bf diagram of topological spaces over} $P$ is
 a~functor $\cd:P\ra\ctopm$.
\end{df}

We shall now consider a~diagram over the face poset $\cf(\da_g)$, to
this end we define a~functor $\cdgn:\cf(\da_g)\ra\ctop$.  Let
$(G,\pi)$ be a~filtered by forests stable graph of genus $g$, and let
$\sigma\in\da_g$ be the~simplex represented by $(G,\pi)$. We set
\begin{equation} \label{eq:cdgn}
\cdgn(\sigma):=\da(G)^n/\aut(G,\pi),
\end{equation}
where $\da(G)^n=\da(G)\times\dots\times\da(G)$ denotes the $n$-fold
direct product of the topological spaces $\da(G)$, and the group
action is the diagonal one:
\[g:(x_1,\dots,x_n)\mapsto (gx_1,\dots,gx_n),\]
for all $x_1,\dots,x_n\in\da(G)$, $g\in\aut(G,\pi)$.  Clearly, the
points of $\da(G)^n/\aut(G,\pi)$ encode all possible marking functions
$p_G:[n]\ra\da(G)$ up to the action of the automorphism group, with
the topology being precisely the one with which we have earlier
endowed the space of marked graphs.

Furthermore, let $\ti\sigma$ be a~simplex on the boundary of $\sigma$,
such that $\dim\ti\sigma+1=\dim\sigma$. This means that $\sigma$
covers $\ti\sigma$ in the poset $\cf(\da_g)$. To define the
corresponding map $\cdgn(\sigma)\ra\cdgn(\ti\sigma)$ we need to
consider two cases. To fix the notations, assume that
$\pi=(E_1,\dots,E_m)$.

\vspace{5pt}

\nin {\bf Case 1.} The simplex $\ti\sigma$ is represented by
a~filtered by forests stable graph $(G,\tau)$, where the ordered
partition $\tau$ is obtained from $\pi$ by merging the blocks $E_i$
and $E_{i+1}$, for some $i=1,\dots,m-1$. In this case $\aut(G,\pi)$ is
a~subgroup of $\aut(G,\tau)$, therefore the induced quotient map
\begin{equation} \label{eq:q1}
q:\da(G)^n/\aut(G,\pi)\ra\da(G)^n/\aut(G,\tau), 
\end{equation}
\[
q:[\,\underline{x}\,]_{\aut(G,\pi)}\mapsto [\,\underline{x}\,]_{\aut(G,\tau)},
\]
is well-defined. Here $\underline{x}=(x_1,\dots,x_n)$, and
$[\,\underline{x}\,]_\Gamma$ denotes the orbit of $\underline{x}$ with
respect to the $\Gamma$-action. We now simply set
$\cdgn(\sigma>\ti\sigma):\cdgn(\sigma)\ra\cdgn(\ti\sigma)$ to be the
map~$q$.

\vspace{5pt}

\nin {\bf Case 2.} The simplex $\ti\sigma$ is represented by
a~filtered by forests stable graph $(H,\tau)$, where the graph $H$ is
obtained from the graph $G$ by shrinking the edges from $E_1$, and the
$\tau=(E_2,\dots,E_m)$ is the corresponding ordered partition. Let
$\psi:\da(G)\ra\da(H)$ denote the induced quotient map, and let
$\psi^n:\da(G)^n\ra\da(H)^n$ its $n$-fold direct product, defined by
$\psi^n(x_1,\dots,x_n):=(\psi(x_1),\dots,\psi(x_n))$. Let $\da(E_1)$
denote the topological union of the edges from $E_1$, we have
$\da(G/E_1)=\da(G)/\da(E_1)$. Since the subspace
$\da(E_1)\subseteq\da(G)$ is invariant under the action of
$\aut(G,\pi)$, every automorphism of $(G,\pi)$ induces an~automorphism
of $(H,\tau)$. We let $\iota:\aut(G,\pi)\hra\aut(H,\tau)$ denote the
corresponding map, which is actually a~group homomorphism. We are now
in a~position to define
\begin{equation} \label{eq:q2}
q:\da(G)^n/\aut(G,\pi)\ra\da(H)^n/\aut(H,\tau), 
\end{equation}
\[
q:[\,\underline{x}\,]_{\aut(G,\pi)}\mapsto [\psi^n(\underline{x})]_{\aut(H,\tau)},
\]
where we use the same notations as in \eqref{eq:q1}. The map $q$ is
well-defined, since for all $g\in\aut(G,\pi)$, and all
$\underline{x}\in\da(G)^n$, we have
\begin{equation}\label{eq:psix}
\psi^n(g(\underline{x}))=\iota(g)(\psi^n(\underline{x})). 
\end{equation}
Also in this case we set
$\cdgn(\sigma>\ti\sigma):\cdgn(\sigma)\ra\cdgn(\ti\sigma)$ to be the
map~$q$.

\subsection{The homotopy colimit presentation} $\,$

\begin{df} \label{df:hocolimdf}
Given a~diagram $\cd$ of topological spaces over a~poset $P$, the {\bf
  homotopy colimit} of $\cd$, denoted $\hocolimm\cd$, is the quotient
space
\[\hocolimm\cd=\coprod_{\sigma=(v_0>\dots> v_n)} 
(\sigma\times\cd(v_0))/\thicksim,\] where the disjoint union is taken
over all chains in $P$. The equivalence relation $\thicksim$ is
generated by: for $\tau_i=(v_0>\dots>\hat v_i>\dots>v_n)$, considered
as a~simplex of $\da(P)$, let $f_i:\tau_i\hra\sigma$ be the inclusion
map, then
\begin{itemize}
\item for $i>0$, $\tau_i\times\cd(v_0)$ is identified with the subset
of $\sigma\times\cd(v_0)$, by the map induced by~$f_i$;
\item for $\tau_0=(v_1>\dots>v_n)$, we have 
$f_0(\alpha)\times x\sim\alpha\times\cd(v_0>v_1)(x)$, for any
$\alpha\in\tau_0$, and $x\in\cd(v_0)$.
\end{itemize}

\end{df}

Given two diagrams $\cd_1$ and $\cd_2$ over the same poset $P$, {\it
  a~diagram map} $\cf$ is a~collection of maps
$\cf(x):\cd_1(x)\ra\cd_2(x)$, for all $x\in P$, which commute with the
diagram structure maps, i.e., for all $x,y\in P$, such that $x>y$, we
have \[\cd_2(x>y)(\cf(x))=\cf(y)(\cd_1(x>y)(x)).\] It is a standard
fact, that a~diagram map induces a~continuous map between the
corresponding homotopy colimits $\cf:\hocolim\cd_1\ra\hocolim\cd_2$.

In particular, taking $\cd_2(x)$ to be a~point, for all $x\in P$, we
get a~diagram of topological spaces, whose colimit is the order
complex $\da(P)$. Setting $\cf(x)$ to be the map which takes
everything to one point, we certainly get a~diagram map. Thus, we
arrive at a~map $p:\hocolim\cd_1\ra\hocolim\cd_2\cong\da(P)$. In these
circumstances, the space $\da(P)$ is called the {\it base space}, the
map $p$ is called the {\it base projection map}, and the preimages of
points under $p$ are the {\it fibers}.  We refer the reader to
\cite[Subsection 15.2.2]{CAT} for further details.

In our specific case we see that the base space is
$X_{g,0}\cong\da_g\cong\bd\da_g\cong\da(\cf(\da_g))$.  It turns out
that the homotopy colimit of the diagram $\cdgn$ describes precisely
the tropical moduli space which we are studying, and that the base
projection map is the natural one.

\begin{thm} \label{thm:main}
The space $\xgn$ is homeomorphic to the homotopy colimit of the
diagram $\cdgn$, whereas the map $\lambda(\emptyset\hra
S):\xgn\ra\da_g$ is the base projection map. In particular, the moduli
space of the rational tropical curves with $n$ marked points $\tmgn$
is homotopy equivalent to $\hocolimm\cdgn$.
\end{thm}
\pr This is pretty much straightforward from our construction. For
a~point $x\in\hocolim\cdgn$, its image $p(x)$ under the base
projection map is represented by a~filtered by forests stable metric
graph $(G,l_G)$, whereas the position of $x$ within its fiber
$\da(G)^n/\aut(G,\pi)$ encodes all the ways to mark $n$ points, modulo
the action of the automorphism group. 

This gives a~map $\rho:\hocolim\cdgn\ra\xgn$. Clearly, that map is
bijective. To see that both $\rho$ and its inverse are continuous we
need to see that small perturbations of a point $x$ in
$\hocolim\cdgn$, resp.\ in $\xgn$, causes small perturbations of
$\rho(x)$, resp.\ of $\rho^{-1}(x)$. There are 3 possibilities for
a~small perturbation of $x\in\hocolim\cdgn$:
\begin{enumerate}
\item[(1)] we stay in the fiber, i.e., the base metric graph does not
  change;
\item[(2)] the base metric graph changes, but we stay within the same
  simplex in $\da_g$;
\item[(3)] we move to an adjacent simplex of higher dimension in $\da_g$.
\end{enumerate}
In $\xgn$ these perturbations correspond to the following.
Perturbation (1) corresponds the graph being fixed and the marks
moving on that graph. In this case the topology is the same as we
mentioned before when discussing the fibers. Perturbation (2)
corresponds to rescaling of the blocks of edges of equal lengths,
together with moving of the marks on the rescaled graph. That topology
is the same in this case as well, follows from our proof of the
homeomorphism $\da_g\cong X_{0,g}$.

When analyzing perturbation (3) we distinguish two cases. In the first
case, no new edges are added, but the blocks of edges of equal lengths
are getting split. Here the topology is the same due to construction
of the attachment map in \eqref{eq:q1}. Finally, in the second case,
we add a~number of very short edges, by replacing some of the vertices
with trees. In that situation the topology is the same due to
construction of the attachment map in \eqref{eq:q2}.  \qed 

\begin{rem}
With the benefit of hindsight, one may now interpret
Theorem~\ref{thm:hmain} and its proof, using the language of homotopy
colimits. This is because already the space $\tmgnb$ can be
represented as a~result of a~gluing construction, similar to homotopy
colimit, whose base would be not $\da_g$, but rather a~difference
between two other generalized simplicial complexes $\da\sm\da'$. The
strong deformation retraction in the proof of Theorem~\ref{thm:hmain}
can then be thought of as corresponding to the retraction of all not
closed cells in $\da\sm\da'$, resulting in the homotopy colimit with
the base~$\da_g$.
\end{rem}

The last theorem of this section describes the standard way the spaces
which are presented as homotopy colimits can be simplified, while
preserving their homotopy type. We shall that technique in
Section~\ref{sect:5}.

\begin{thm} \label{thm:homlem} {\rm (Homotopy Lemma)}

\nin Let $\cf:\cd_1\ra\cd_2$ be a~diagram map between diagrams of
spaces over~$\da$, such that for each $v\in\da^{(0)}$, the map
$\cf(v):\cd_1(v)\ra\cd_2(v)$ is a~homotopy equivalence. Then the
induced map $\hocolimm\cf:\hocolimm\cd_1\ra\hocolimm\cd_2$ is
a~homotopy equivalence as well.
\end{thm}

Again, we refer to~\cite[Chapter 15]{CAT} for the proof and further
information.

\newpage

\section{CW structure on $\xgn$ and its Euler characteristic}
\label{sect:cw} \label{sect:4}

We shall now describe CW structure on the spaces $\xgn$ derived from
the homotopy colimit representation from Theorem~\ref{thm:main}.  

\subsection{Cubical complexes associated to filtered stable graphs} $\,$

\nin To start with, we associate cubical complexes to all simplices
of~$\da_g$. For a~filtered stable graph $(G,\pi)$ we let $S(G,\pi)$
denote the $1$-dimensional simplicial complex obtained from $\da(G)$
by inserting the middle point into every edge of $G$ which is flipped
by some element of the group $\aut(G,\pi)$; this of course should
include all the loops of~$G$.

\begin{figure}[hbt]
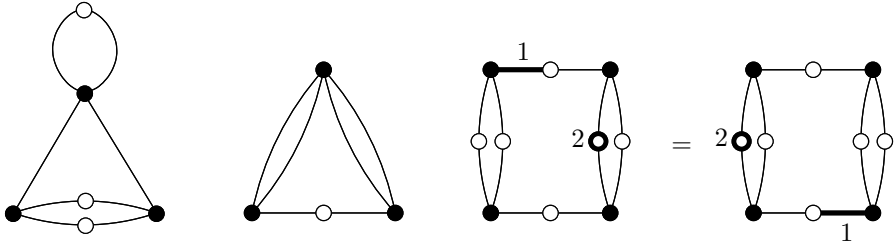

\begin{center}
  \begin{picture}(0,0)%
    \includegraphics{f5.pstex}%
  \end{picture}%
  \input{f5.pstex_t}%
  
\end{center}
\caption{Here we show three examples of the subdivisions $S(G,\pi)$,
  in each one we have $\pi=(E(G))$. We depict the original vertices of
  $G$ as filled-in and the added subdivision points as hollow. On the
  right hand side we show examples of two different presentations of
  the same cube in $C(G,\pi)$; here the fattened-up half-edge and
  a~vertex indicate the allowed loci of the marked points. }
\label{fig:f5}
\end{figure}

\begin{df} \label{df:cgpi}
Let $\sigma$ be a~simplex of $\da_g$ represented by a~filtered stable
graph $(G,\pi)$ of genus~$g$. We set
\begin{equation} \label{eq:cgpi}
C(G,\pi):=S(G,\pi)^n/\aut(G,\pi),
\end{equation}
where the action of the group $\aut(G,\pi)$ on the direct product
$S(G,\pi)^n$ is the diagonal one.
\end{df}
We shall also use the notation $C(\sigma)$. Since $S(G,\pi)$ is
a~subdivision of $\da(G)$, the comparison of the definitions
\eqref{eq:cdgn} and \eqref{eq:cgpi} shows that $C(G,\pi)$ is
homeomorphic to $\cdgn(G,\pi)$. We shall think of each cube of
$S(G,\pi)^n$ as specifying where every label $i=1,\dots,n$ should lie.
We shall call the corresponding part of $S(G,\pi)$ the {\it allowed
  locus of $i$}, and note that it is either a~vertex of $S(G,\pi)$, or
an edge or a~half-edge of~$G$. The same is true for $C(G,\pi)$, but of
course modulo the $\aut(G,\pi)$-action. See Figure~\ref{fig:f5} for
examples.

Importantly, the space $C(G,\pi)$ has the structure of a~cubical
complex. This is because $S(G,\pi)$ is a~cubical complex, and the
$\aut(G,\pi)$-action on $S(G,\pi)^n$ has the following property: if
the cube of $S(G,\pi)^n$ is preserved by an~element $g\in\aut(G,\pi)$,
then it must be fixed by it pointwise. Indeed, if $g\tau=\tau$, for
some $g\in\aut(G,\pi)$, and a~cube $\tau$ from $S(G,\pi)^n$, then for
all $i=1,\dots,n$, the allowed locus of $i$ is preserved
by~$g$. Clearly, if the allowed locus of $i$ is a~vertex of
$S(G,\pi)$, or a~half-edge of $G$, then it must be preserved
pointwise. If, on the other hand, the allowed locus is an~edge $e$ of
$G$ and it is not preserved pointwise by $g$, then $g$ must flip $e$,
hence, by our construction, the edge $e$ should have been subdivided
in $S(G,\pi)$, yielding a~contradiction.


\subsection{The CW structure on $\xgn$} $\,$


\nin{\it The indexing of cells.} The cells are indexed by pairs
$(\sigma,\tau)$, where $\sigma$ is a~simplex of the generalized
simplicial complex $\da_g$, and $\tau$ is a~cube of $C(\sigma)$. We
shall denote such a~cell $c(\sigma,\tau)$. The open cell
$c(\sigma,\tau)$ consists of all the points of $\xgn$, which are
indexed by filtered by forests stable graphs $(G,\pi)$ representing
the simplex $\sigma$ of $\da_g$, with $n$ marked points, such that the
point marked~$i$ belongs to the corresponding allowed locus prescribed
by the cube~$\tau$, for $i=1,\dots,n$, modulo the action of the group
$\aut(G,\pi)$. The dimension of the cell $c(\sigma,\tau)$ is equal to
$\dim\sigma+\dim\tau$. 

For example, the cell depicted on the right hand side of
Figure~\ref{fig:f5} is indexed by the pair $(\sigma,\tau)$, where
$\sigma$ is a~vertex, indexed by $(G,\pi)$, where $G$ is the graph
with $4$ vertices and $6$ edges shown there, and $\pi=(E(G))$, and
$\tau$ is a~$1$-dimensional cube. In particular, this cell has
dimension~$1$. Of course, the two presentations of this cell may
correspond to different cells, if the ordered partition $\pi$ is
chosen differently.

\vspace{5pt}

\nin{\it The attachment maps.} Consider the cell
$c(\sigma,\tau)$. Assume the simplex $\sigma$ is represented by
a~filtered by forests stable graph $(G,\pi)$, where
$\pi=(E_1,\dots,E_t)$. To describe the attachment of $c(\sigma,\tau)$
we need to tell how to glue in the cells from the boundary $\bo
c(\sigma,\tau)$; for this it is enough to take the cells of dimension
$\dim\sigma+\dim\tau-1$ from $(\bo\sigma)\times\tau$ and
$\sigma\times(\bo\tau)$. The case when we take the cells from
$\sigma\times(\bo\tau)$ is easy, as we simply glue along the
attachment map of the cell $\tau$ in the cubical complex $C(\sigma)$.

Let us consider the case when we take a~cell $\lambda\times\tau$ from
$(\bo\sigma)\times\tau$ of dimension $\dim\sigma+\dim\tau-1$. There
are two distinguished cases.

\vspace{5pt}

\noindent {\bf Case 1.} The simplex $\lambda$ is indexed by the
filtered by forests stable graphs $(G,\pi')$ obtained from $(G,\pi)$
by merging blocks $E_i$ and $E_{i+1}$, for some $i=1,\dots,t$,
in~$\pi$. Then, we have $\aut(G,\pi)\subseteq\aut(G,\pi')$, which
induces the quotient map
\[q:S(G,\pi)^n/\aut(G,\pi)\lra S(G,\pi')^n/\aut(G,\pi'),\]
defined by $q:[x]\mapsto[x]$, where we identify the spaces $S(G,\pi)$
and $S(G,\pi')$. This is precisely the gluing map from that part of
boundary of $c(\sigma,\tau)$ to $c(\lambda,\tau)$ which we are looking
for.

\vspace{5pt}

\noindent {\bf Case 2.} The simplex $\lambda$ is indexed by the
filtered by forests stable graphs $(G',\pi')$ obtained from $(G,\pi)$
by shrinking the edges $E_1$. In that case, we have a~shrinking map
$\psi^n:S(G,\pi)^n\ra S(G',\pi')^n$, and the map
$\iota:\aut(G,\pi)\hra\aut(G',\pi')$, which together induce the
quotient cubical map 
\[q:S(G,\pi)^n/\aut(G,\pi)\lra S(G',\pi')^n/\aut(G',\pi'),
\]
defined by $q:[x]\mapsto[\psi^n(x)]$. This map is well-defined due to
the identity~\eqref{eq:psix}, and provides us with the required gluing
map from that part of boundary of $c(\sigma,\tau)$
to~$c(\lambda,\tau)$.

\vspace{5pt}

This cell structure allows us to write down the chain complex
$(C_*(\xgn;\zz);\bo)$ of vector spaces over $\zz$ for computing
$H_*(\xgn;\zz)$. The vector space $C_t(\xgn;\zz)$ is generated by all
the $t$-cells, and the algebraic boundary maps are induced by the
topological boundary described above, with every individual gluing
contributing the cell one is glued onto with coefficient $1$, if the
cell has the right dimension. The degenerate gluings yield
contribution~$0$. We remark that the gluing is degenerate if and only
if we are gluing a~cell of dimension $\dim\sigma+\dim\tau-1$ from
$(\bo\sigma)\times\tau$, we are in Case 2 above, and some of the $n$
marked points happen to lie on one of the edges from the set $E_1$.

\subsection{The colimit of $\cdgn$} $\,$

\nin It is curious, that while the homotopy colimit of the diagram
$\cdgn$ is a~rather complicated space, its colimit is in fact quite
simple.

\begin{prop}
The colimit of the diagram $\cdgn$ is a~point.
\end{prop}
\pr One can describe the space $\colim\cdgn$ explicitly as follows:
it is the union of all the spaces in the diagram $\cdgn$ modulo
a~certain equivalence relation. This relation is generated by the
elementary relations, which say that for
$\alpha,\alpha'\in\cf(\da_g)$, such that $\alpha>\alpha'$, and
$x\in\cdgn(\alpha)$, we have $x\sim\cdgn(\alpha>\alpha')(x)$. 

Let us now fix $\alpha$ to be the vertex of $\da_g$ indexed by the
graph with one vertex and $g$ loops, and let $w$ denote that unique
vertex. The $n$-tuple $\bar w=(w,\dots,w)$ indexes a~vertex in the
cubical complex $C(\alpha)$, and hence also a~point in $\colim\cdgn$,
we shall show that every other point
$x\in\cup_{\sigma\in\da_g}\cdgn(\sigma)$ is equivalent to $\bar w$.
Using the described cubical structures on $C(\sigma)$, it is actually
enough to show that any cube in any $C(\sigma)$ is equivalent to
$\bar w$.

Let $c$ be a cube of some $C(G,\pi)$, and assume that $c$ is indexed
by the $n$-tuple $(a_1,\dots,a_n)$, where each $a_1$ is either
a~vertex or an~edge of $S(G,\pi)$. To start with, setting all edges of
$G$ equal corresponds to a~map in the diagram $\cdgn$, and of course
the cube $c$ will be equivalent to its image under this map.
Therefore, we may replace $c$ with its image, and assume that all the
edges of $G$ have the same length.  

Assume now that some $a_i$ is either a~barycenter of a~loop $l$ of
$\da(G)$, or a~half-edge lying on a~loop $l$ of $G$. In that case,
there exists a~filtered by forests stable graph $(G',\pi')$, such that
shrinking the shortest edges of $G'$ will yield $(G,\pi)$. For
example, it can be obtained by any admissible ``unlooping'' of $l$.
That shrinking corresponds to a~surjective map in the diagram $\cdgn$,
and we can replace the cube $c$ with any of its preimage cubes in
$C(G',\pi')$, which is of course equivalent to $c$. Repeating this
procedure, we can see that we may assume that no $a_i$ is a~part of
a~loop of $G$.

Assume now, that some $a_i$ is an~edge of $S(G,\pi)$, and let $e$
denote the underlying edge of $G$. Let $(G',\pi')$ be obtained from
$(G,\pi)$ by letting $e$ be slightly shorter than the other
edges. There is cube $c'\in C(G',\pi')$, indexed by
$(a_1',\dots,a_n')$ which is mapped to $c$ by the corresponding map
$f$ of $\cdgn$, hence $c\sim c'$. Let $c''$ be a~cube of $C(G',\pi')$
obtained from $c'$ be replacing $a_i'$ by any of its endpoints.  We
can shrink $e$ in $G'$ obtaining $(G'',\pi'')$. Under the map of the
diagram $\cdgn$ corresponding to this shrinking, the cubes $c'$ and
$c''$ map to the same cube, hence they are equivalent. We conclude
that $c$ is equivalent to $f(c'')$.  Repeating this procedure, we can
see that we may assume that all $a_i$'s correspond to vertices of
$S(G,\pi)$.

Clearly, the argument of the last paragraph can be also used to show
that $c$ is equivalent to $c'=(a_1',\dots,a_n')$, where the $n$-tuple
$(a_1',\dots,a_n')$ is obtained from $(a_1,\dots,a_n)$ be replacing
any $a_i$ by any neighboring vertex of $S(G,\pi)$. Repeating this, we
end up with the cube where all the $a_i$'a are the same vertices of
$S(G,\pi)$ an correspond to a~vertex of $G$. Taking the appropriate
shrinking map to $\alpha$ we conclude that all the cubes are
equivalent to the vertex $\bar w$, and hence the entire colimit is
just a~point.  \qed

\subsection{General formula for the generating function of the numbers
of cells of $\xgn$} $\,$

\nin We shall now use the cell structure of $\xgn$ to derive formulae
for its Euler characteristic. For this we will need the following
weighted version of Burnside Lemma, see e.g., \cite[Exercise 14.4.5,
  page 313]{Biggs}.

\begin{lm} \label{lm:wbl} {\rm (Weighted Burnside Lemma). }

\nin Let $\Gamma$ be a group acting on a set $X$, and let
$w:X\ra\Omega$ be a~$\Gamma$-invariant map, where $(\Omega,+)$ is
an~abelian group, then we have
\begin{equation} \label{eq:wbl}
\sum_{o\in O(\Gamma,X)}w(o)=
\frac{1}{|\Gamma|}\sum_{g\in\Gamma}\sum_{x\in\fixm(g)}w(x),
\end{equation}
where $O(\Gamma,X)$ denotes the set of orbits of the $\Gamma$-action
on $X$.
\end{lm}

For a~finite-dimensional CW complex $K$ and a~variable $x$, we let
$\cp(K)(x)$ denote the generating function for the numbers of cells,
i.e., $\cp(K):=\sum_{d=0}^{\dim K}n_d x^d$, where $n_d$ denotes the
number of $d$-cells in $K$. When the choice of the variable is clear,
we shall simply write $\cp(K)$. Note that $\cp(K)(-1)$ equals to the
(nonreduced) Euler characteristic of~$K$.

Let now $X$ be the set of the cubes of $S(G,\pi)^n$, and let $\Gamma$
be the group $\aut(G,\pi)$ with the standard $\Gamma$-action on
$X$. Let us furthermore set $w(\sigma):=x^{\dim\sigma}$ for all cubes
$\sigma$ of $S(G,\pi)^n$. Since $\sum_{d=0}^{\dim K}n_d
x^d=\sum_{\sigma\in K}x^{\dim\sigma}$, for any CW complex $K$, we get
$\cp(C(G,\pi))=\sum_{\sigma\in C(G,\pi)}w(\sigma)$. Substituting this
data into~\eqref{eq:wbl} we obtain
\begin{equation} \label{eq:wbl2}
\cp(C(G,\pi))=
\frac{1}{|\aut(G,\pi)|}\sum_{\gamma\in\aut(G,\pi)}\cp(\fix_{C(G,\pi)}(\gamma)),
\end{equation}
where $\fix_{C(G,\pi)}(\gamma)$ denotes the subcomplex of $C(G,\pi)$
fixed by the group element~$\gamma$. It follows from our construction
that $\fix_{C(G,\pi)}(\gamma)=(\fix_{S(G,\pi)}(\gamma))^n$, where
$\fix_{S(G,\pi)}(\gamma)$ denotes the subcomplex of $S(G,\pi)$ fixed
by the group element~$\gamma$. Since $\cp(K^n)=\cp(K)^n$ for an
arbitrary finite CW complex $K$, the equation~\eqref{eq:wbl2}
translates to equation
\begin{equation} \label{eq:wbl3}
\cp(C(G,\pi))=
\frac{1}{|\aut(G,\pi)|}\sum_{\gamma\in\aut(G,\pi)}
\cp(\fix_{S(G,\pi)}(\gamma))^n.
\end{equation}
In particular, substituting $x=-1$ into \eqref{eq:wbl3} we obtain the
equation for the corresponding Euler characteristic
\begin{equation} \label{eq:wbl4}
\chi(C(G,\pi))=
\frac{1}{|\aut(G,\pi)|}\sum_{\gamma\in\aut(G,\pi)}
\chi(\fix_{S(G,\pi)}(\gamma))^n.
\end{equation}

Summing with appropriate signs over all the isomorphism classes of
filtered by forests stable graphs, equations \eqref{eq:wbl3} and
\eqref{eq:wbl4} yield equations
\begin{multline} \label{eq:wbl5}
\cp(\xgn)=\sum_{\sigma\in\da_g}(-1)^{\dim\sigma}\cp(C(\sigma))= \\
=\frac{1}{|\aut(G,\pi)|}\sum_{\sigma\in\da_g}(-1)^{\dim\sigma}
\sum_{\gamma\in\aut(G,\pi)}\cp(\fix_{S(\sigma)}(\gamma))^n.
\end{multline}
and
\begin{equation} \label{eq:wbl6}
\chi(\xgn)=
\frac{1}{|\aut(G,\pi)|}\sum_{\sigma\in\da_g}(-1)^{\dim\sigma}
\sum_{\gamma\in\aut(G,\pi)}\chi(\fix_{S(G,\pi)}(\gamma))^n.
\end{equation}

These formulae will come in handy in our analysis of cases of small
genus.

\section{The case of genus $2$} \label{sect:5}

In this and the next section we will analyse the cases of small genus
$g=2$ and $g=3$, and make some conjectures which should hold for all
$g$ and $n$. 

\subsection{The cases $n=0$ and $n=1$} $\,$

\begin{figure}[hbt]
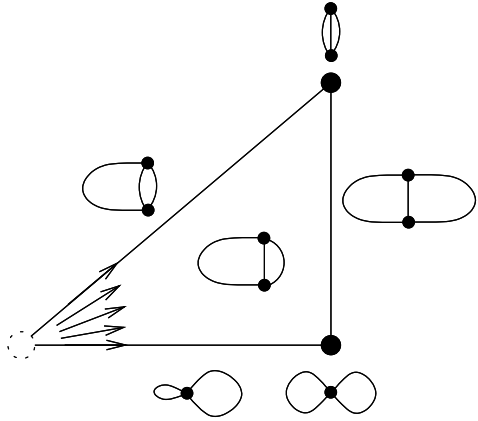

\begin{center}
  \begin{picture}(0,0)%
    \includegraphics{g21.pstex}%
  \end{picture}%
  \input{g21.pstex_t}%
  
\end{center}
\caption{The strong deformation retraction of $TM_{2,n}^b$ onto $TM_{2,n}^e$
as seen in the projection forgetting the marked points.}
\label{fig:g21}
\end{figure}

\nin Let us start with considering the space $TM_{2,n}$. As was shown
in~\cite{tms}, this space is homotopy equivalent to
$TM_{2,n}^b$. Furthermore, by Theorem~\ref{thm:hmain}, there is
a~strong deformation retraction from $TM_{2,n}^b$ to $TM_{2,n}^b$; it
is symbolically shown on Figure~\ref{fig:g21} with arrows.

\begin{figure}[hbt]
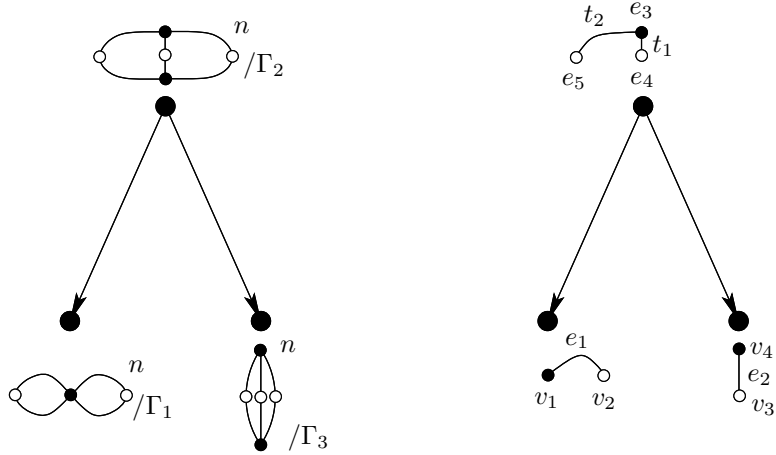

\begin{center}
  \begin{picture}(0,0)%
    \includegraphics{g22.pstex}%
  \end{picture}%
  \input{g22.pstex_t}%
  
\end{center}
\caption{On the left hand side we show the diagram $\cd_{2,n}$, and
  the right hand side we show the special case $n=1$.}
\label{fig:g22}
\end{figure}

As in general case, it is now sufficient to understand the space
$X_{2,n}$.  According to Theorem~\ref{thm:main} that space is
a~homotopy colimit shown on the left hand side of
Figure~\ref{fig:g22}. The groups on that figure are as follows:
\begin{itemize}
\item the group $\Gamma_1$ consists of the flips of the loops, with
  a~possible swap of the two loops; it has cardinality $8$ and is
  isomorphic to the wreath product of $\cs_2$ with $\cs_2$;
\item the group $\Gamma_2$ is generated by reflections of the graph
  with respect to the horizontal and vertical axes; it has cardinality
  $4$ and is isomorphic to the direct product of $\cs_2$ with $\cs_2$;
\item the group $\Gamma_3$ is generated by $\cs_3$ permuting the edges
  and the reflection with respect to the horizontal axis; it has
  cardinality $12$ and is isomorphic to the direct product of $\cs_3$
  with $\cs_2$.
\end{itemize}
We again used the hollow circles to denote the points inserted in the
midpoints of edges, as is prescribed by the group actions; in fact
here, we simply need to subdivide all the edges. As elaborated on in
Section~\ref{sect:cw} these subdivisions induce a~CW structure on the
homotopy colimit. For $n=0$, this structure coincides with the
generalized simplicial complex $\da_g$, so here we simply get
a~$1$-simplex. 

\begin{figure}[hbt]
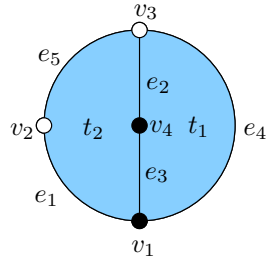

\begin{center}
  \begin{picture}(0,0)%
    \includegraphics{g23.pstex}%
  \end{picture}%
  \input{g23.pstex_t}%
  
\end{center}
\caption{The CW structure on $X_{2,1}$.}
\label{fig:g23}
\end{figure}

For $n=1$ we get the diagram of spaces shown on the right hand side of
Figure~\ref{fig:g22}. In that diagram, the map corresponding to the
diagonal arrow pointing southwest shrinks the edge $t_1$ to the point
$v_1$, and maps the edge $t_2$ homeomorphically to the edge $e_1$,
whereas the map corresponding to the diagonal arrow pointing southeast
takes the vertex $e_3$ to $v_4$, and maps both edges $t_1$ and $t_2$
homeomorphically to the edge $e_2$. The homotopy colimit of this
diagram is the cell complex of dimension $2$, shown on
Figure~\ref{fig:g23}. Clearly, that complex is collapsible in the
sense of \cite{cohen}, in particular, it is contractible.

\subsection{The case $n=2$} $\,$

\nin For shorthand notations, let $A\stackrel f\leftarrow C\stackrel
g\rightarrow B$ denote the spaces and maps in the diagram on the left
hand side of Figure~\ref{fig:g22}, for $n=2$, with the space $A$
corresponding the space in the southwest corner of that diagram. We
start by understanding the cubical structure of the spaces $A$, $B$,
and $C$, all of which are $2$-dimensional pure cubical complexes.

\begin{figure}[hbt]
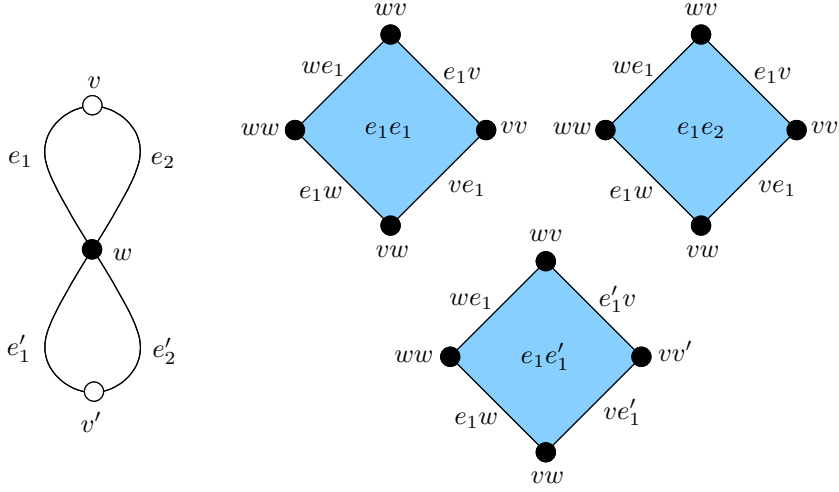

\begin{center}
  \begin{picture}(0,0)%
    \includegraphics{g2n2a.pstex}%
  \end{picture}%
  \input{g2n2a.pstex_t}%
  
\end{center}
\caption{The cubical complex $A$.}
\label{fig:g2n2a}
\end{figure}

The case of the cubical complex $A$ is illustrated on
Figure~\ref{fig:g2n2a}. Here the filtered by forests stable graph
$(G,\pi)$ consists of one vertex and two loops of length $1$. The
subdivided space $S(G,\pi)$ is shown on the left hand side of
Figure~\ref{fig:g2n2a}, where also notations for vertices and edges
are fixed. Using these notations we describe the cubes of $A$. Each
one is described by two letters, corresponding to the allowed loci of
the points marked $1$ and $2$ (in that order). Thus, the complex $A$
has

\nin \underline{$5$ vertices:} $ww$, $wv$, $vw$, $vv$, $vv'$;

\nin \underline{$6$ edges:} $we_1$, $e_1w$, $ve_1$, $e_1v$, $ve_1'$, $e_1'v$; 

\nin \underline{$3$ $2$-cubes:} $e_1e_1$, $e_1e_2$, $e_1e_1'$;

\nin where we of course only list one of the equivalent descriptions
for each cube, for example $vv'=v'v$, and $e_1'v=e_1v'=e_2v'=e_2'v$.

We see from the Figure~\ref{fig:g2n2a} that the $2$-cubes $e_1 e_1$ and
$e_1 e_2$ have common boundary, and hence form a~$2$-sphere; and that
the $2$-cube $e_1 e_1'$ is attached to that $2$-sphere along two
neighboring sides. We conclude that $A$ is homotopy equivalent to
a~$2$-sphere.

\begin{figure}[hbt]
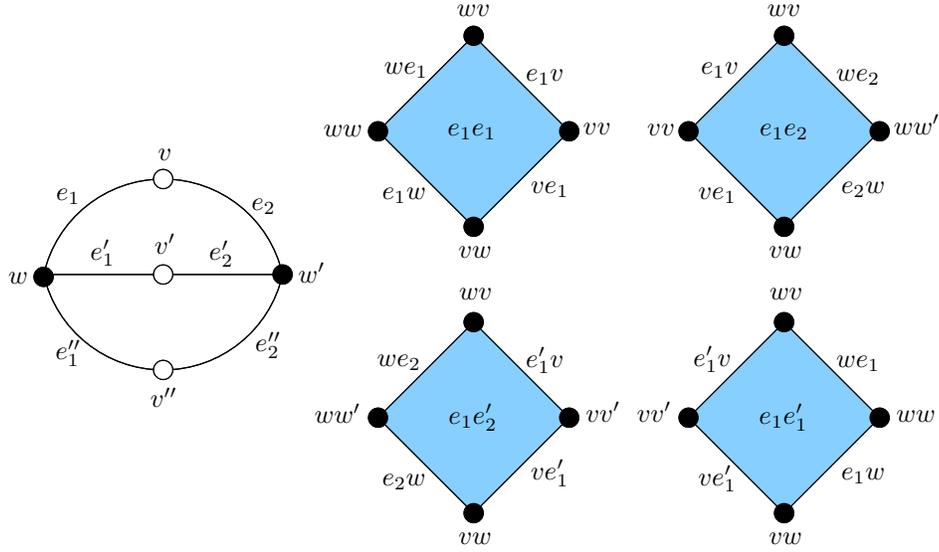

\begin{center}
  \begin{picture}(0,0)%
    \includegraphics{g2n2b.pstex}%
  \end{picture}%
  \input{g2n2b.pstex_t}%
  
\end{center}
\caption{The cubical complex $B$.}
\label{fig:g2n2b}
\end{figure}

Next we describe the cubical complex $B$. The corresponding subdivided
filtered by forests stable graph is shown on the left hand side of
Figure~\ref{fig:g2n2b}. With the notations there, the complex $B$ has

\nin \underline{$6$ vertices:} $ww$, $ww',$ $wv$, $vw$, $vv$, $vv'$;

\nin \underline{$8$ edges:} $we_1$, $e_1w$, $we_2$, $e_2w$, 
$ve_1$, $e_1v$, $ve_1'$, $e_1'v$;

\nin \underline{$4$ $2$-cubes:} $e_1e_1$, $e_1e_2$, $e_1e_1'$, $e_1e_2'$;

\nin where we again list only one of the equivalent descriptions
for each cube. The $4$ $2$-cubes can be glued to each other in the order
on the Figure~\ref{fig:g2n2b} forming a~cubical complex homeomorphic
to the $2$-sphere. That complex may also conveniently be visualised as
follows: take the double pyramid with a $2$-cube base (also known as
a~$3$-dimensional crosspolytope, which is the polytope dual to the
$3$-cube) and delete the $4$ edges of the base, each time merging the
adjacent triangles to form $2$-cubes.

\begin{figure}[hbt]
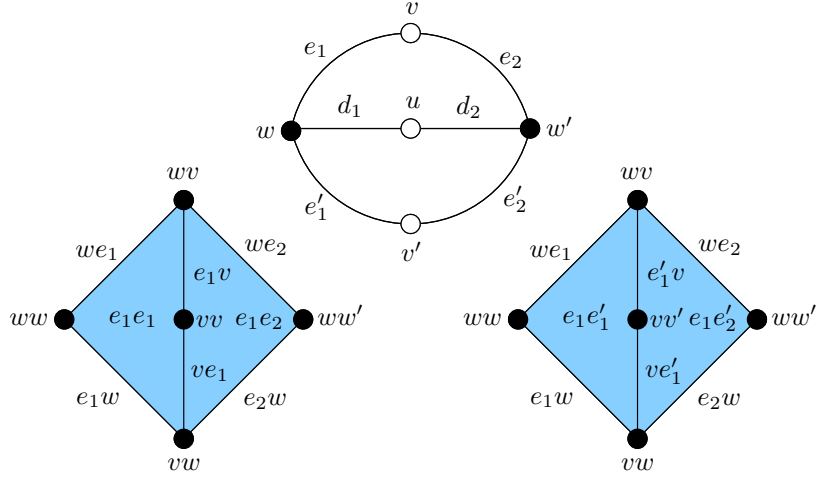

\begin{center}
  \begin{picture}(0,0)%
    \includegraphics{g2n2c.pstex}%
  \end{picture}%
  \input{g2n2c.pstex_t}%
  
\end{center}
\caption{The subdivided filtered by forests stable graph corresponding
  to the cubical complex $C$, together with a~part of that complex.}
\label{fig:g2n2c}
\end{figure}

Finally, we describe the cubical complex $C$. The corresponding
subdivided filtered by forests stable graph is shown in the middle of
Figure~\ref{fig:g2n2c}. With the notations there, the complex $C$ has

\nin \underline{$11$ vertices:} $uu$, $uw$, $uv$, $wu$, $vu$,
$ww$, $ww',$ $wv$, $vw$, $vv$, $vv'$;

\nin \underline{$18$ edges:} $ud_1$, $d_1u$, $ue_1$, $e_1u$, 
$vd_1$, $d_1v$, $wd_1$, $d_1w$, $wd_2$, $d_2w$, $we_1$, $e_1w$,
$we_2$, $e_2w$, $ve_1$, $e_1v$, $ve_1'$, $e_1'v$;

\nin \underline{$10$ $2$-cubes:} $d_1d_1$, $d_1e_1$, $d_1e_2$, 
$d_1d_2$, $e_1d_1$, $e_1d_2$, $e_1e_1$, $e_1e_2$, $e_1e_1'$,
$e_1e_2'$.

\begin{figure}[hbt]
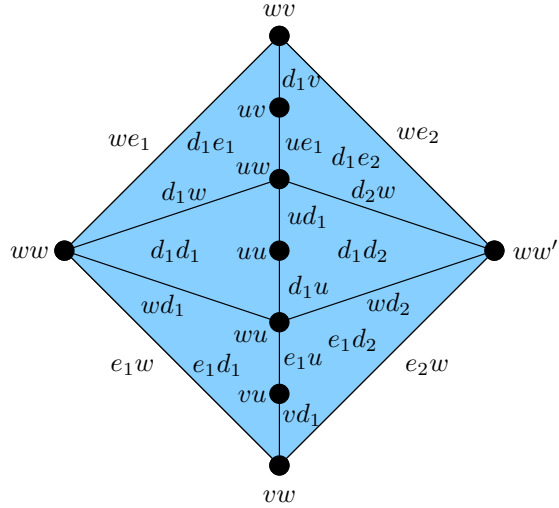

\begin{center}
  \begin{picture}(0,0)%
    \includegraphics{g2n2c2.pstex}%
  \end{picture}%
  \input{g2n2c2.pstex_t}%
  
\end{center}
\caption{Six $2$-cells glued to form a~single $2$-cell inside the
  cubical complex $C$.}
\label{fig:g2n2c2}
\end{figure}

As Figures~\ref{fig:g2n2c} and~\ref{fig:g2n2c2} illustrate, it is
helpful to divide the $10$ constituting $2$-cubes into $3$ groups:
$e_1e_1$ and $e_1e_2$, $e_1e_1'$ and $e_1e_2'$, and the remaining $6$
$2$-cubes. The $2$-cubes in each group, when glued, form a~$2$-cube,
and all these $3$ new $2$-cubes have a~common boundary. Thus, we
conclude that the complex $C$ is homeomorphic to a~$2$-sphere with
a~disc glued in, which in turn, is homotopy equivalent to a~wedge of
two $2$-spheres $\cs^2\vee\cs^2$.

As part of the analysis of the maps $f:C\ra A$ and $g:C\ra B$, we show
in Table~\ref{tab:1} what their values on $2$-cubes; there we use $*$
to denote lower-dimensional cubes. We invite the reader to get the
geometric intuition for what these maps actually do to the spaces.

\begin{table} [hbt]
\[
\begin{array} {c|c|c}
\text{the $2$-cube of $C$} & \text{image under $f$} & 
\text{image under $g$} \\
\hline
e_1e_1 & e_1e_1 & e_1e_1 \\
e_1e_2 & e_1e_2 & e_1e_2 \\
e_1e_1' & e_1e_1' & e_1e_1' \\
e_1e_2' & e_1e_1' & e_1e_2' \\
d_1d_1 & * & e_1e_1 \\
d_1e_1 & * & e_1e_1' \\
d_1e_2 & * & e_1e_2' \\
d_1d_2 & * & e_1e_2 \\
e_1d_1 & * & e_1e_1' \\
e_1d_2 & * & e_1e_2' 
\end{array}
\]
\caption{The values of $f$ and $g$ on the $2$-cubes of $C$.}
\label{tab:1}
\end{table}

We are now ready to prove the following result.

\begin{thm}
The moduli space of tropical curves of genus $2$ with two marked
points $TM_{2,2}$ is contractible.
\end{thm}

\pr Let us describe subcomplexes of $A$, $B$, and $C$:
\begin{itemize}
\item the complex $\wti A$ is obtained from $A$ by deleting the
  interior of the $2$-cell $e_1e_1$;
\item the complex $\wti B$ is obtained from $A$ by deleting the
  interior of the $2$-cell $e_1e_1$;
\item the complex $\wti C$ is obtained from $A$ by deleting the
  interiors of the $2$-cells $e_1e_1$ and $d_1d_1$.
\end{itemize}
It follows from our detailed descriptions of the cubical structures on
$A$, $B$, and $C$, that we have inclusions $f(\wti C)\subseteq\wti A$
and $g(\wti C)\subseteq\wti B$. It follows from the standard
properties of maps of quotients, see e.g., \cite[Satz 11.7,
  p.\ 115]{rinow}, that the induced maps $\ti f:C/\wti C\ra A/\wti A$
and $\ti g:C/\wti C\ra B/\wti B$ are well-defined by $\ti
f([x]):=[f(x)]$ and $\ti g([x]):=[g(x)]$. Furthermore, the quotient
maps $q^A:A\ra A/\wti A$, $q^B:B\ra B/\wti B$, and $q^C:C\ra C/\wti
C$, induce a~diagram map $q$ from $A\stackrel f\la C\stackrel g\ra B$
to $A/\wti A\stackrel{\ti f}\la C/\wti C\stackrel{\ti g}\ra B/\wti B$.

As is clear from our explicit description, all three complexes $\wti
A$, $\wti B$, and $\wti C$, are contractible. Since inclusions of CW
subcomplexes are cofibrations, see e.g., \cite[Proposition 0.17]{hat},
it follows that the maps $q^A$, $q^B$, and $q^C$, are homotopy
equivalences. Using Theorem~\ref{thm:homlem} (Homotopy Lemma), we
conclude that $q$ induces a~homotopy equivalence between the homotopy
colimits of these two diagrams.

\begin{figure}[hbt]
\begin{center}
  \begin{picture}(0,0)%
    \includegraphics{g2n2f.pstex}%
  \end{picture}%
  \input{g2n2f.pstex_t}%
  
\end{center}
\caption{On the left hand side we show the diagram $A/\wti
  A\stackrel{\ti f}\la C/\wti C\stackrel{\ti g}\ra B/\wti B$, whereas
  on the right hand side we show how its homotopy colimit can be
  obtained by a~simple self-identification on a~cone over~$\cs^2$.}
\label{fig:g2n2f}
\end{figure}

By our construction and previous developments the homotopy colimit of
$A\stackrel f\la C\stackrel g\ra B$ is homotopy equivalent to the
tropical module space $TM_{2,2}$. On the other hand, the spaces
$A/\wti A$ and $B/\wti B$ are each homeomorphic to a~$2$-sphere,
whereas the space $C/\wti C$ is homeomorphic to a~wedge of $2$-spheres
$\cs^2\vee\cs^2$. The map $\ti f$ shrinks one of the spheres in the
wedge to a~point and maps the other one homeomorphically to
$\cs^2\cong A/\wti A$, while the map $\ti g$ maps both $2$-spheres in
$C$ homeomorphically to $\cs^2\cong B/\wti B$. Hence the homotopy
colimit of $A/\wti A\stackrel{\ti f}\la C/\wti C\stackrel{\ti g}\ra
B/\wti B$ is homeomorphic to the space obtained from a~cone over
$\cs^2$ by taking an edge connecting the apex to one of the points in
the base, dividing this edge in the middle and identifying the two
halves with each other by a~reflection about the middle point, see
Figure~\ref{fig:g2n2f}. In particular, it is contractible, and we have
therefore proved the theorem. \qed

\subsection{The Euler characteristic of the space $X_{2,n}$ and 
the proof that this space is almost never contractible} $\,$

\begin{figure}[hbt]
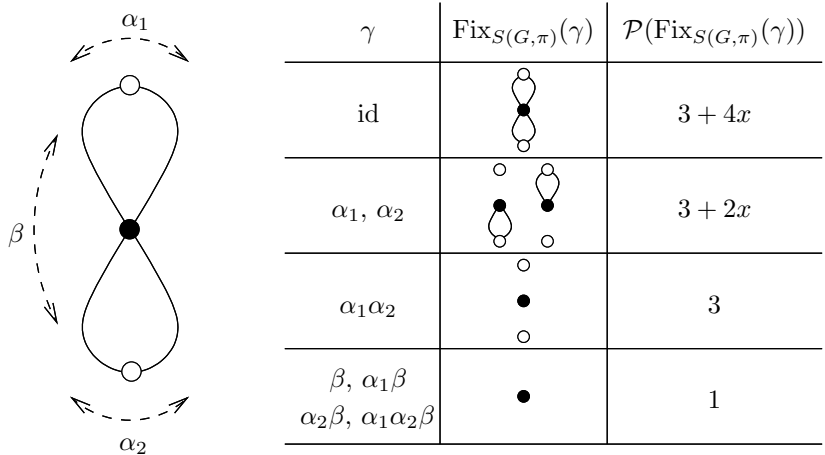

\begin{center}
  \begin{picture}(0,0)%
    \includegraphics{mf1.pstex}%
  \end{picture}%
  \input{mf1.pstex_t}%
  
\end{center}
\caption{Both loops of the graph on the left hand side have length~$1$. 
The elements $\alpha_1$ and $\alpha_2$ are reflections of the
corresponding loops, and $\beta$ is swapping the loops.}
\label{fig:mf1}
\end{figure}

\begin{figure}[hbt]
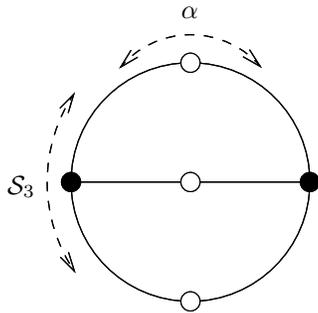

\begin{center}
  \begin{picture}(0,0)%
    \includegraphics{mf2.pstex}%
  \end{picture}%
  \input{mf2.pstex_t}%
  
\end{center}
\caption{The uppermost and the lowest edge of the graph on the left hand
side have equal length, while the middle edge is shorter. The element
$\alpha$ corresponds to the vertical reflection, whereas the element
$\beta$ swaps the two long edges.}
\label{fig:mf2}
\end{figure}

\nin Let us now calculate of the Euler characteristic, and more generally,
the generating function for the number of cells of $X_{2,n}$. This can
be done using the formula~\eqref{eq:wbl5}. For this, we need to
calculate $\cp(\fix_{S(G,\pi)}(\gamma))$, for all filtered by forests
stable graphs of genus $2$, and all $\gamma\in\aut(G,\pi)$. That
information is summarized in tables on Figures~\ref{fig:mf1},
\ref{fig:mf2}, and \ref{fig:mf3}, corresponding to the three different
filtered by forests stable graphs of genus~$2$. On
Figure~\ref{fig:mf3} we use the convention that
\[0^n=\begin{cases}
1,& \text{ for } n=0,\\
0,& \text{ for } n\geq 1.
\end{cases}
\]

\begin{figure}[hbt]
\begin{center}
  \begin{picture}(0,0)%
    \includegraphics{mf3.pstex}%
  \end{picture}%
  \input{mf3.pstex_t}%
  
\end{center}
\caption{In the graph on the left hand side all the three edges are of 
length~$1$. Again, the element $\alpha$ corresponds to the vertical
reflection. In addition, we have a~permutation action of $\cs_3$ on
the edges.}
\label{fig:mf3}
\end{figure}

We conclude that 
\begin{multline}\label{eq:chi3}
\cp(X_{2,n})=\frac{1}{24}(3\cdot(3+4x)^n-4\cdot(5+6x)^n+6\cdot(3+2x)^n
-3^n+\\ +4\cdot2^n+12+4\cdot 0^n),
\end{multline}
and hence
\begin{equation}\label{eq:chi2}
\chi(X_{2,n})=\frac{1}{24}\left(-3^n+(-1)^{n+1}+2^{n+2}+18+4\cdot 0^n\right),
\end{equation}
for all nonnegative integers $n$. For small values of $n$ we get
\[\chi(X_{2,0})=\chi(X_{2,1})=\chi(X_{2,2})=\chi(X_{2,3})=1,
\quad\chi(X_{2,4})=0,\quad\chi(X_{2,5})=-4.\] In particular, we
conclude that $X_{2,n}$ is not contractible for all $n\geq 4$.  At
present, it is unknown whether $X_{2,3}$ is contractible or not.

\newpage $\,$ \newpage

\section{The case of genus $3$} \label{sect:g3} \label{sect:6}

\subsection{Collapsibility of $\da_3$} $\,$

\begin{figure}[hbt]
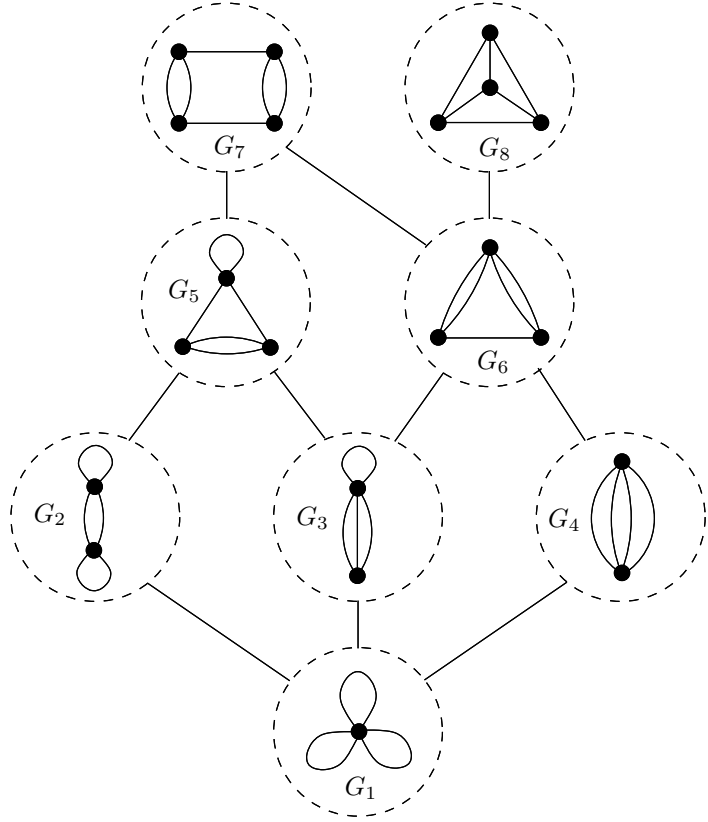

\begin{center}
  \begin{picture}(0,0)%
    \includegraphics{g31.pstex}%
  \end{picture}%
  \input{g31.pstex_t}%
  
\end{center}
\caption{The eight isomorphism classes of stable graphs of genus~$3$,
  corresponding to the eight vertices of~$\da_3$.}
\label{fig:g31}
\end{figure}

\nin As was shown in Proposition~\ref{prop:3.4}, the generalized
simplicial complex $\da_3$ is pure, and it has dimension~$4$.  There
are eight isomorphism classes of stable graphs of genus~$3$. These are
shown on Figure~\ref{fig:g31}, and correspond to eight vertices
of~$\da_3$.

\begin{thm}
The generalized simplicial complex $\da_3\cong X_{3,0}$ is
collapsible.
\end{thm}

\pr As a~first simplification, we notice that the vertex $G_2$ belongs
to a~unique $3$-simplex. This $3$-simplex is indexed by the filtered
by forests stable graph $(G_7,\pi)$, with
$\pi=(\{e_1\},\{e_2\},\{e_3\},E(G_7)\sm\{e_1,e_2,e_3\})$, where $e_1$
is one of the vertical edges, $e_2$ is another vertical edge,
non-adjacent to $e_1$, and $e_3$ is a~horizontal edge connecting the
two. Therefore, deletion of the vertex $G_2$ along with all the
adjacent simplices corresponds to a~(non-elementary) simplicial
collapse. Let us denote by $X$ the thus obtained generalized
simplicial complex. It is enough to show that this complex is
collapsible.

\begin{figure}[hbt]
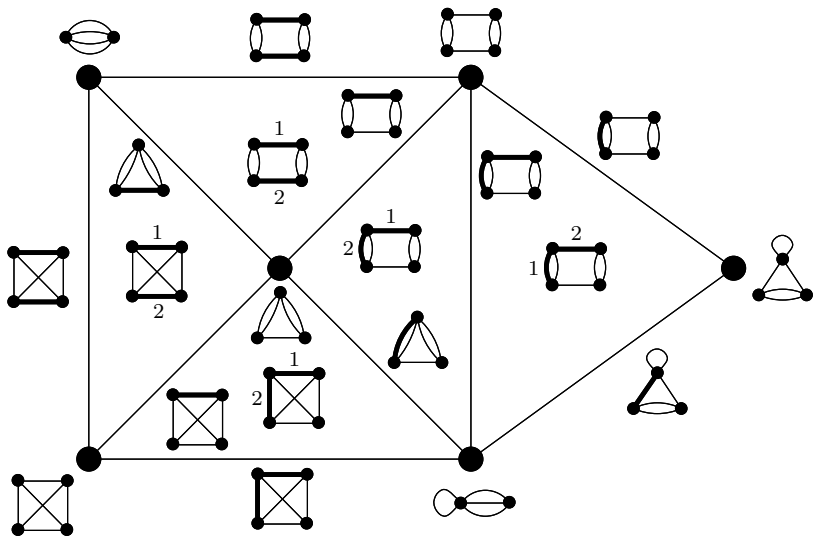

\begin{center}
  \begin{picture}(0,0)%
    \includegraphics{g32.pstex}%
  \end{picture}%
  \input{g32.pstex_t}%
  
\end{center}
\caption{The deletion of the vertex $G_1$ from the complex $X$.}
\label{fig:g32}
\end{figure}

\begin{figure}[hbt]
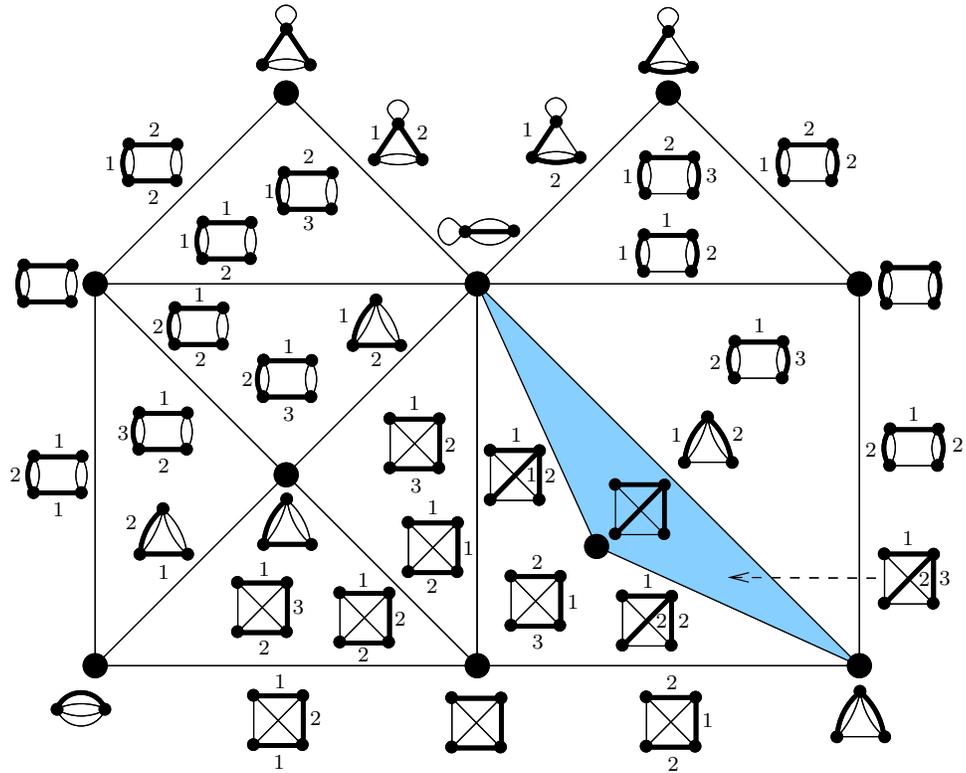

\begin{center}
  \begin{picture}(0,0)%
    \includegraphics{g33.pstex}%
  \end{picture}%
  \input{g33.pstex_t}%
  
\end{center}
\caption{The link of the vertex $G_1$ in the complex $X$.}
\label{fig:g33}
\end{figure}

We notice that the closed star of the vertex $G_1$ constitutes the
entire complex $X$. Hence, it is enough to show that the generalized
simplicial complexes $\dl_X(G_1)$, the deletion of $G_1$ from $X$, and
$\lk_X(G_1)$, the link of $G_1$ in $X$, are both collapsible. These
complexes are shown on Figures~\ref{fig:g32} and \ref{fig:g33}, from
which it is apparent that there are various ways to collapse each one
to a~point. \qed

\subsection{The asymptotics of Euler characteristic of $X_{3,n}$, and 
the conjectural asymptotics for all $\xgn$} $\,$

\begin{figure}[hbt]
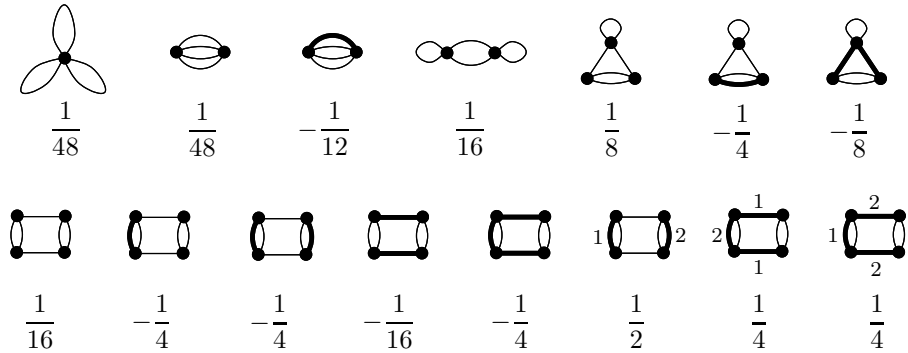

\begin{center}
  \begin{picture}(0,0)%
    \includegraphics{g34.pstex}%
  \end{picture}%
  \input{g34.pstex_t}%
  
\end{center}
\caption{The contributions to the coefficient of $4^n$ to the Euler
  characteristic of $TM_{3,n}$ made by the Euler characteristics of
  the various cells of $\da_3$.}
\label{fig:g34}
\end{figure}

\nin Let us now show that $\Omega(\chi(X_{3,n}))=4^n$, that is
$\Omega(\chi(X_{3,n}))/4^n$ converges to a~constant, as $n$ goes to
infinity. We use the formula~\eqref{eq:wbl6}. The nontrivial
contributions from different cells of $\da_3$ are summarized on
Figure~\ref{fig:g34}, where we list the contributions below
corresponding cells. Summing up we see that the final answer
is~$1/48$.

It was shown in \cite{tms2} that for $n\geq 1$, we have
$TM_{1,{n+1}}\simeq(\cs^1)^n/\zz$, where the $\zz$-action is
a~simultaneous reflection on each factor circle $\cs^1$. Taking the
points fixed by that reflection as vertices, and two semicircles as
edges we induce a~$\zz$-invariant cubical structure on $(\cs^1)^n$
with $\cp((\cs^1)^n)=(2+2x)^n$. Using equation~\eqref{eq:wbl} we
obtain the equality $\cp((\cs^1)^n/\zz)=\frac{1}{2}((2+2x)^n+2^n)$. In
particular, we have $\chi(TM_{1,n})=2^{n-2}$, for $n\geq 2$.

Emboldened by the settled cases $g=1,2$, and $3$, we arrive at the
following conjecture.
\begin{conj} \label{conj:as}
Assume that $g$ be a~fixed arbitrary positive integer. We have
\begin{equation}\label{eq:asy}
\Omega(\chi(TM_{g,n}))=(g+1)^n,
\end{equation}
where $\chi(TM_{g,n})$ is considered as a~function of $n$, and the
asymptotics is taken for $n\ra\infty$.
\end{conj}

It would also be interesting to find an~interpretation for the actual
coefficient of $(g+1)^n$, the values which we have computed are
$1/4$, $-1/24$, $1/48$.

\section{Conclusion and open questions}

Though having gained a~substantial advance in our understanding of the
structure of the spaces $\xgn$, and thus of the homotopy type of the
moduli spaces of rational tropical curves $\tmgn$, quite a~few
questions remain open.

\begin{figure}[hbt]
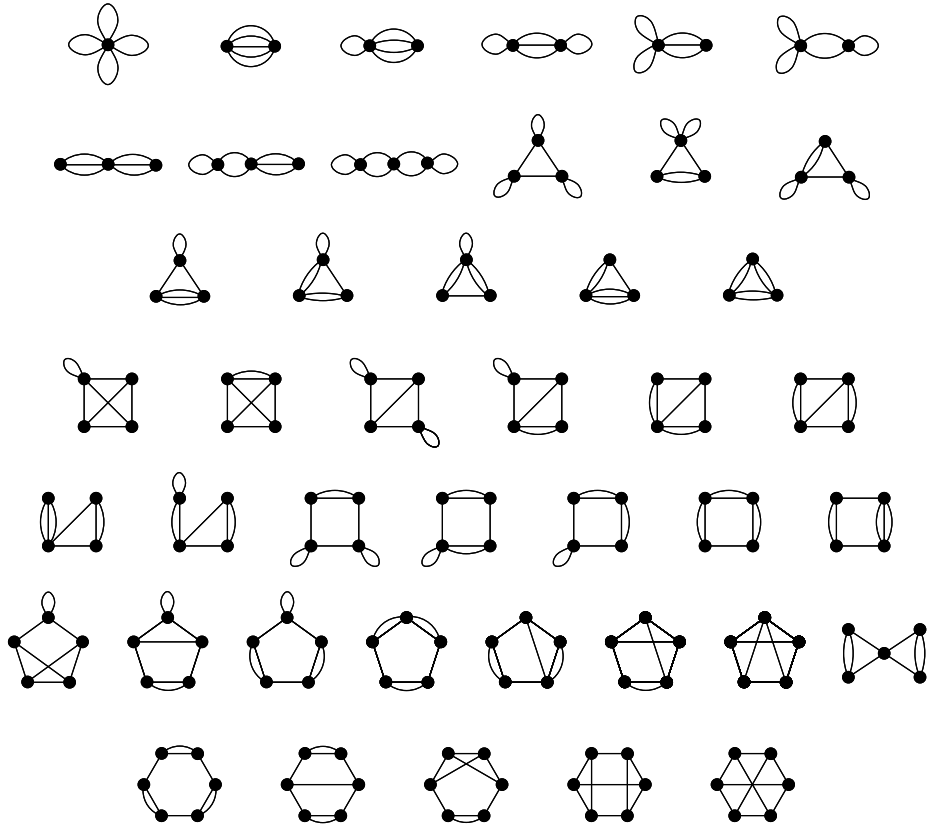

\begin{center}
  \begin{picture}(0,0)%
    \includegraphics{g4.pstex}%
  \end{picture}%
  \input{g4.pstex_t}%
  
\end{center}
\caption{The $43$ isomorphism classes of stable graphs of genus~$4$,
  corresponding to the $43$ vertices of~$\da_4$.}
\label{fig:g4}
\end{figure}

For example, it seems natural to conjecture that $\da_g$ is not always
contractible, and it would be interesting to find out whether that
phenomenon occurs already for $g=4$. At the moment, not much is known
about $\da_4$ beyond the fact that it has $43$ vertices, whose
indexing stable graphs of genus $4$ are shown on Figure~\ref{fig:g4}.

Furthermore, adding the marked points into the fray, makes for a~lot
more questions. As a~first concrete one, we would like to know whether
$X_{2,3}$ is contractible or not.  In general, it certainly feels like
the spaces $\xgn$ should always never be contractible, so settling
Conjecture~\ref{conj:as} seems useful.

In Section~\ref{sect:5} we could understand the homotopy type of the
space $X_{2,2}$ by coming up with some deformations of the spaces in
the diagram $\cdgn$, and using the Homotopy Lemma. The natural
question here is of course whether that technique could be turned into
a~system and help us understand other instances of $g$ and~$n$.

\vskip5pt

\nin {\bf Acknowledgments.} The author would like to thank Eva-Maria
Feichtner for the valuable discussions. He would also like to thank
the University of Bremen for the support of the work on this project
within the framework of AG CALTOP. Part of this research has been done
while the author was visiting the Mathematical Research Institute at
Oberwolfach, which is hereby gratefully acknowledged.

\end{document}